\documentclass[12pt]{article}
\usepackage[english,french]{babel}
\usepackage{latexsym}
\usepackage{amsmath}
\usepackage{amssymb}
\usepackage{graphics}
\usepackage{graphicx}
\usepackage{fancyhdr}
\usepackage{color}
\usepackage{latexsym}
\usepackage{epstopdf}
\usepackage{amsfonts}
\usepackage[utf8]{inputenc}
\usepackage[numbers]{natbib}				
\renewcommand{\cite}[1]{\citeauthor*{#1} [\citeyear{#1}]}
\usepackage{amsthm,url,xspace}
\definecolor{couleurCitations}{rgb}{0,0.65,0}
\definecolor{couleurRef}{rgb}{0.75,0,0}
\usepackage[pdftex,,nativepdf=true,colorlinks=true,
linkcolor=couleurRef,citecolor=blue,a4paper=true
,bookmarks=true,plainpages=true,
 urlcolor= blue] {hyperref}

\title{\textbf{Pointwise Adaptive M-estimation in Nonparametric Regression}}
\author{\textbf{Michaël CHICHIGNOUD$^1$}}
\date{}
\makeatletter

\@addtoreset{equation}{section}
\makeatother

\begin{document}

\theoremstyle{plain}
\newtheorem{theorem}{Theorem}
\newtheorem{proposition}{Proposition}
\newtheorem{lemma}{Lemma}
\newtheorem{definition}{Definition}
\newtheorem{remark}{Remark}
\newtheorem{assumptions}{Assumption}
\newtheorem{properties}{Properties}
\newtheorem{corollary}{Corollary}



\def\restriction#1#2{\mathchoice
              {\setbox1\hbox{${\displaystyle #1}_{\scriptstyle #2}$}
              \restrictionaux{#1}{#2}}
              {\setbox1\hbox{${\textstyle #1}_{\scriptstyle #2}$}
              \restrictionaux{#1}{#2}}
              {\setbox1\hbox{${\scriptstyle #1}_{\scriptscriptstyle #2}$}
              \restrictionaux{#1}{#2}}
              {\setbox1\hbox{${\scriptscriptstyle #1}_{\scriptscriptstyle #2}$}
              \restrictionaux{#1}{#2}}}
\def\restrictionaux#1#2{{#1\,\smash{\vrule height .8\ht1 depth .85\dp1}}_{\,#2}}

\newcommand{\epr}{\hfill\hbox{\hskip 4pt
                \vrule width 5pt height 6pt depth 1.5pt}\vspace{0.5cm}\par}

\newcommand{\rd}{{\rm d}}
\newcommand{\rk}{\textsf{k}}
\newcommand{\rf}{\textsf{f}}
\newcommand{\rh}{\textsf{h}}
\newcommand{\rz}{\textsf{z}}
\newcommand{\rK}{\textsf{K}}
\newcommand{\rQ}{\textsf{Q}}
\newcommand{\rc}{\textsf{c}}
%
%
\newcommand{\al}{\alpha}
\newcommand{\bt}{\beta}
\newcommand{\g} {\gamma}
\newcommand{\dl}{\delta}
\newcommand{\sg}{\sigma}
\newcommand{\ze}{\zeta}
\newcommand{\Dl}{\Delta}
\newcommand{\te}{\theta}
\newcommand{\vth}{\vartheta}
\newcommand{\vf}{\varphi}
\newcommand{\ups}{\upsilon}
\newcommand{\Up}{\Upsilon}
\newcommand{\ve}{\varepsilon}
\newcommand{\e}{\varepsilon}
%
%
\newcommand{\cA}{{\cal A}}
\newcommand{\cB}{{\cal B}}
\newcommand{\cC}{{\cal C}}
\newcommand{\cD}{{\cal D}}
\newcommand{\cE}{{\cal E}}
\newcommand{\cF}{{\cal F}}
\newcommand{\cG}{{\cal G}}
\newcommand{\cH}{{\cal H}}
\newcommand{\cI}{{\cal I}}
\newcommand{\cJ}{{\cal J}}
\newcommand{\cK}{{\cal K}}
\newcommand{\cL}{{\cal L}}
\newcommand{\cM}{{\cal M}}
\newcommand{\cN}{{\cal N}}
\newcommand{\cO}{{\cal O}}
\newcommand{\cP}{{\cal P}}
\newcommand{\cQ}{{\cal Q}}
\newcommand{\cR}{{\cal R}}
\newcommand{\cS}{{\cal S}}
\newcommand{\cT}{{\cal T}}
\newcommand{\cU}{{\cal U}}
\newcommand{\cV}{{\cal V}}
\newcommand{\cW}{{\cal W}}
\newcommand{\cX}{{\cal X}}
\newcommand{\cY}{{\cal Y}}
\newcommand{\cZ}{{\cal Z}}
\newcommand{\cz}{{\cal z}}
\newcommand{\calr}{{\cal r}}
%
\newcommand{\Ba}{\boldsymbol{a}}
\newcommand{\BA}{\boldsymbol{A}}
\newcommand{\Bb}{\boldsymbol{b}}
\newcommand{\BB}{\boldsymbol{B}}
\newcommand{\Bc}{\boldsymbol{c}}
\newcommand{\BC}{\boldsymbol{C}}
\newcommand{\Bd}{\boldsymbol{d}}
\newcommand{\BD}{\boldsymbol{D}}
\newcommand{\Be}{\boldsymbol{e}}
\newcommand{\BE}{\boldsymbol{E}}
\newcommand{\Bf}{\boldsymbol{f}}
\newcommand{\BF}{\boldsymbol{F}}
\newcommand{\Bg}{\boldsymbol{g}}
\newcommand{\BG}{\boldsymbol{G}}
\newcommand{\Bh}{\boldsymbol{h}}
\newcommand{\BH}{\boldsymbol{H}}
\newcommand{\Bi}{\boldsymbol{i}}
\newcommand{\BI}{\boldsymbol{I}}
\newcommand{\Bj}{\boldsymbol{j}}
\newcommand{\BJ}{\boldsymbol{J}}
\newcommand{\Bk}{\boldsymbol{k}}
\newcommand{\BK}{\boldsymbol{K}}
\newcommand{\Bl}{\boldsymbol{l}}
\newcommand{\BL}{\boldsymbol{L}}
\newcommand{\Bm}{\boldsymbol{m}}
\newcommand{\BM}{\boldsymbol{M}}
\newcommand{\Bn}{\boldsymbol{n}}
\newcommand{\BN}{\boldsymbol{N}}
\newcommand{\Bo}{\boldsymbol{o}}
\newcommand{\BO}{\boldsymbol{O}}
\newcommand{\Bp}{\boldsymbol{p}}
\newcommand{\BP}{\boldsymbol{P}}
\newcommand{\Bq}{\boldsymbol{q}}
\newcommand{\BQ}{\boldsymbol{Q}}
\newcommand{\Br}{\boldsymbol{r}}
\newcommand{\BR}{\boldsymbol{R}}
\newcommand{\Bs}{\boldsymbol{s}}
\newcommand{\BS}{\boldsymbol{S}}
\newcommand{\Bt}{\boldsymbol{t}}
\newcommand{\BT}{\boldsymbol{T}}
\newcommand{\Bu}{\boldsymbol{u}}
\newcommand{\BU}{\boldsymbol{U}}
\newcommand{\Bv}{\boldsymbol{v}}
\newcommand{\BV}{\boldsymbol{V}}
\newcommand{\Bw}{\boldsymbol{w}}
\newcommand{\BW}{\boldsymbol{W}}
\newcommand{\Bx}{\boldsymbol{x}}
\newcommand{\BX}{\boldsymbol{X}}
\newcommand{\By}{\boldsymbol{y}}
\newcommand{\BY}{\boldsymbol{Y}}
\newcommand{\Bz}{\boldsymbol{z}}
\newcommand{\BZ}{\boldsymbol{Z}}
%
%
\newcommand{\Balpha}{\boldsymbol{\alpha}}
\newcommand{\Bbeta}{\boldsymbol{\beta}}
\newcommand{\Bgamma}{\boldsymbol{\gamma}}
\newcommand{\Bdelta}{\boldsymbol{\delta}}
\newcommand{\Btheta}{\boldsymbol{\theta}}
\newcommand{\BTheta}{\boldsymbol{\Theta}}
\newcommand{\Btau}{\boldsymbol{\tau}}
\newcommand{\Bla}{\boldsymbol{\lambda}}
\newcommand{\Bmu}{\boldsymbol{\mu}}
\newcommand{\Blal}{\boldsymbol{\lambda_1}}
\newcommand{\Blad}{\boldsymbol{\lambda_2}}
\newcommand{\Bxi}{\boldsymbol{\xi}}
%
%
\newcommand{\bA}{\mathbb A}
\newcommand{\bB}{\mathbb B}
\newcommand{\bC}{\mathbb C}
\newcommand{\bD}{\mathbb D}
\newcommand{\bE}{\mathbb E}
\newcommand{\bF}{\mathbb F}
\newcommand{\bG}{\mathbb G}
\newcommand{\bH}{\mathbb H}
\newcommand{\bI}{{\mathbb I}}
\newcommand{\bJ}{{\mathbb J}}
\newcommand{\bK}{{\mathbb K}}
\newcommand{\bL}{{\mathbb L}}
\newcommand{\bM}{{\mathbb M}}
\newcommand{\bN}{{\mathbb N}}
\newcommand{\bO}{{\mathbb O}}
\newcommand{\bP}{{\mathbb P}}
\newcommand{\bQ}{{\mathbb Q}}
\newcommand{\bR}{{\mathbb R}}
\newcommand{\bS}{{\mathbb S}}
\newcommand{\bT}{{\mathbb T}}
\newcommand{\bU}{{\mathbb U}}
\newcommand{\bV}{{\mathbb V}}
\newcommand{\bW}{{\mathbb W}}
\newcommand{\bX}{{\mathbb X}}
\newcommand{\bY}{{\mathbb Y}}
\newcommand{\bZ}{{\mathbb Z}}
%
%
\newcommand{\mA}{\mathfrak{A}}
\newcommand{\mB}{\mathfrak{B}}
\newcommand{\mF}{\mathfrak{F}}
\newcommand{\mH}{\mathfrak{H}}
\newcommand{\mS}{\mathfrak{S}}
\newcommand{\mL}{\mathfrak{L}}
\newcommand{\ms}{\mathfrak{s}}
\newcommand{\mP}{\mathfrak{P}}
\newcommand{\mK}{\mathfrak{K}}
\newcommand{\mM}{\mathfrak{M}}
\newcommand{\mk}{\mathfrak{k}}
\newcommand{\mI}{\mathfrak{I}}
\newcommand{\mW}{\mathfrak{W}}
\newcommand{\mi}{\mathfrak{i}}
\newcommand{\mh}{\mathfrak{h}}
\newcommand{\mj}{\mathfrak{j}}
\newcommand{\mw}{\mathfrak{w}}
\newcommand{\mJ}{\mathfrak{J}}
\newcommand{\mm}{\mathfrak{m}}
\newcommand{\mn}{\mathfrak{n}}
\newcommand{\mv}{\mathfrak{v}}
\newcommand{\mT}{\mathfrak{T}}
\newcommand{\mX}{\mathfrak{X}}
\newcommand{\mZ}{\mathfrak{Z}}

%
%

%
%
\newcommand{\estimatorb}[1]{\hat{#1}}			
\newcommand{\estimatorh}[1]{\breve{#1}}			
\newcommand{\estimatorlp}[1]{\mathring{#1}}		
\newcommand{\estimatorml}[1]{\dot{#1}}			
\newcommand{\Ef}[1]{\bE_{#1}} 				
\newcommand{\design}[1]{X_{#1}}				
\newcommand{\observ}[1]{Y_{#1}}				
\newcommand{\noiseUniform}[1]{U_{#1}}			
\newcommand{\density}{g}				
\newcommand{\densityAlpha}{g_\alpha}			
\newcommand{\noiseAdditif}[1]{\xi_{#1}}			
\newcommand{\noiseAlpha}[1]{\epsilon_{#1}}			
\newcommand{\point}{y}					
\newcommand{\minimaxrate}[1]{\varphi_n({#1})}		
\newcommand{\minimaxrateGamma}[1]{\varphi_{n,\gamma}({#1})}
\newcommand{\fenetre}{h}				
\newcommand{\voisin}[2]{V_{#1}(#2)}			
\newcommand{\numberderivateb}{D_b}			
\newcommand{\setcoefficient}{\Theta}			
\newcommand{\likelihood}{L_\fenetre}			
\newcommand{\normalisation}[1]{N_{#1}}			
\newcommand{\adaptiverate}[1]{\phi_n({#1})}		
\newcommand{\pricepay}[1]{\rho_n({#1})}			
\newcommand{\adaptiverateGamma}[1]{\phi_{n,\gamma}({#1})}
\newcommand{\pricepayGamma}[1]{\rho_{n,\gamma}({#1})}	


\def\restriction#1#2{\mathchoice
              {\setbox1\hbox{${\displaystyle #1}_{\scriptstyle #2}$}
              \restrictionaux{#1}{#2}}
              {\setbox1\hbox{${\textstyle #1}_{\scriptstyle #2}$}
              \restrictionaux{#1}{#2}}
              {\setbox1\hbox{${\scriptstyle #1}_{\scriptscriptstyle #2}$}
              \restrictionaux{#1}{#2}}
              {\setbox1\hbox{${\scriptscriptstyle #1}_{\scriptscriptstyle #2}$}
              \restrictionaux{#1}{#2}}}
\def\restrictionaux#1#2{{#1\,\smash{\vrule height .8\ht1 depth .85\dp1}}_{\,#2}}




\maketitle
\selectlanguage{english}
\begin{center}
$^1$\emph{Université Aix-Marseille 1, LATP\\
39, rue F. Joliot Curie, 13453 Marseille Cedex 13, France.\\
E-mail: chichign@cmi.univ-mrs.fr}
\end{center}

\medskip
\begin{abstract}
This paper deals with the nonparametric estimation in heterosce\-dastic
regression $ Y_i=f(X_i)+\xi_i, \: i=1,\ldots,n $, with
incomplete
information, i.e. each real random variable $ \xi_i $ has a
density $ g_{i} $ which is unknown to
the statistician. The aim is to estimate the regression function $ f $ at a given point. Using a local polynomial
fitting from 
M-estimator denoted $ \estimatorb f^h $ and applying Lepski's procedure for the bandwidth
selection, we construct an estimator $ \estimatorb f^{\hat h} $
which is adaptive over the collection of isotropic H\"{o}lder classes. In particular, we establish new exponential
inequalities to
control deviations of local M-estimators allowing to construct the minimax estimator. The advantage of this
estimator is that
it does not depend on densities of
random errors and we only assume
that the probability density functions are symmetric and monotonically on $ \bR_+ $. It is important to mention that our
estimator is robust compared to extreme values of the noise.
\\
\\
{\bf Key words:} adaptation, Huber function, Lepski's method, M-estimation,
minimax estimation, nonparametric regression, robust estimation,
pointwise estimation.
\\
\\
{\bf AMS 2000 subject classification:}{ Primary 62G08; secondary 62G20, 62G35.}
\end{abstract}

\section{Introduction}

Let the statistical experiment be generated by the
observation \\$Z^{(n)}=(X_i,Y_i)_{i=1,...n}, n\in\bN^{*}$,
where each $(X_i,Y_i)$ satisfies the equation
\begin{equation}
\label{model} Y_i=f(X_i)+\xi_i,\quad i=1,\ldots, n.
\end{equation}
Here $f:[0,1]^d\rightarrow\mathbb{R}$ is an unknown function to be estimated at a given point $x_0\in[0,1]^d$ from the
observation $Z^{(n)}$.

The real random variables $(\xi_i)_{i\in{1,\ldots,n}}$ (the noise) are supposed to be independent
and each variable $ \xi_i $ has a symmetric density
$g_i(\cdot) $, with respect to the Lebesgue measure on
$\mathbb R$. We also assumed that $g_i$ is monotonically on $ \bR_+ $ for any $i$.

The design points $(X_i)_{i\in{1,...,n}}$ are independent and uniformly distributed on $[0,1]^d$. The
random vectors
$(X_i)_{i\in{1,...,n}}$ and $(\xi_i)_{i\in{1,\ldots,n}}$ are independent.

Along the paper, the unknown function $f$ is supposed to
be smooth,
in particular, it belongs to an isotropic Hölder ball of functions
$\mathbb{H}_d({\beta},L,M)$ (cf. Definition \ref{def_holder}
below).
Here ${\beta}>0$ is the smoothness of $f$, $L>0$ is the Lipschitz
constant and $M$ is an upper bound of $f$ and its partial derivatives.

\paragraph{Motivation.} In this paper, the considered problem is the robust nonparametric estimation, i.e. the
estimation of the regression function $f$ in the presence of a heavy-tailed noise (cf.
\cite{Rousseuw_Leroy87} and
\cite{Huber_Ronchetti09}). Well-known examples are when the noise distribution is for instance Laplace (no finite
exponential's moment) or Cauchy (no finite order's moments).
Moreover, we assume that the noise
densities $ \big(g_i\big)_{i=1,...,n} $ are unknown to the statistician. This problem has popular applications, for
example in relative GPS positioning (cf.
\cite{Chang_Guo}) or in robust image denoising (cf. \cite{Katkovnik_Foi_Egiazarian_Astola10}).

In parametric case, we consider $ f $ as a constant parameter $ \theta\in\bR $. The use of empiric criteria is
very
popular, i.e. the minimization of the following {\it contrast function} $ \rho $:
$$
\hat\theta=\arg\min_{t\in[-M,M]}\sum_{i=1}^n\rho(Y_i-t),
$$ 
The most famous contrast functions are the square function $ \rho(z)=z^2 $ ($\hat\theta$ become the empiric mean), the
absolute value function $ \rho(z)=|z| $ ($\hat\theta$ become the empiric mean)
and the {\it Huber function}, as defined in (\ref{huber_function}), without an explicit expression of $ \hat\theta $
(cf. \cite{Huber64}).
It is well known that the square function
leads to the empiric mean which does not fit with a heavy-tailed noise. Thus the square function is not suitable in the
model (\ref{model}).

\bigskip

In nonparametric estimation, we propose a {\it local parametric approach} (LPA) to
estimate the regression function at a given point $ x_0\in[0,1]^d $ in the model (\ref{model}). We suppose that $ f $ is
locally almost polynomial
(with degree $b\in\bN$) and we use the parametric estimator on a neighborhood denoted $ V_{x_0}(h) $
The parameter is reconstructed from the following
criterion, for any $x_0\in[0,1]^d$ and $ h\in[0,1]$
\begin{equation}\label{Intro_criterion}
 \hat\theta=\arg\min_{t\in[-M,M]^{N_b}}\sum_{i:X_i\in V_{x_0}(h)
}\rho\big(Y_i-f_t(X_i)\big)\:K\left(\frac{X_i-x_0}{h}\right).
\end{equation}
where $ f_t(\cdot) $ is a polynomial of degree $b$ with coefficients $ t $, $K(\cdot)$ is a {\it kernel function} and $
N_b $ is the number of partial
derivatives of $ f $ of order smaller than $b$.
We refer to $\hat f^h(x_0)=f_{\hat\theta}(x_0)$ as the {\it $\rho$-LPA
estimator}. It belongs to the family
of M-estimators and it relies on a local scale parameter $h$, called the {\it bandwidth}. A crucial issue is the optimal
choice of the parameter $h$. To adress it we use quite standard arguments based on the {\it bias/variance trade-off}
(cf. (\ref{equation_minimiser}) below) in minimax case and the {\it Lepski's rule} for the {\it data-driven selection}
in adaptation. First, since $ f
$ is smooth ($
f\in\mathbb{H}_d({\beta},L,M) $,
cf. Definition \ref{def_holder} below) we notice that
\begin{equation}\label{def_term_bias}
\exists\theta=\theta(f,x_0,h)\in \Theta\big(M\big):\quad b_h:=\sup_{x\in V_{x_0}(h)}\big|f(x)-f_\theta(x)\big|\leq
Ldh^\beta.
\end{equation}
We can choose $ \theta $ as the coefficients of Taylor
polynomial as defined in (\ref{coefficient_taylor}).	
Thus, if $h$ is chosen sufficiently small our original model
(\ref{model}) is well approximated inside of $V_{x_0}(h)$ by the
``parametric" model
\begin{equation}\label{model_parametric}
 \cY_i=f_\theta(X_i)+\xi_i,\quad
\forall i\::\:X_i\in V_{x_0}(h).
\end{equation}
With this model, the {\it $ \rho $-LPA estimator} $\hat{\theta}$
achieves the usual parametric rate of convergence $ 1/\sqrt{nh^d} $, where $ nh^d $ is the number of the observations
in the neighborhood $V_{x_0}(h)$ (See Theorem \ref{minimax}, Section \ref{sectionmaximalrisk}).

This approach has been introduced by \cite{Katkovnik85} and used for the first time in robust nonparametric estimation
by \cite{Tsybakov86},
\cite{Hardle_Tsybakov88} and \cite{Hall_Jones89} to obtain asymptotic normality and minimax results. We also notice
that \citeauthor{Tsybakov82a}[\citeyear{Tsybakov82a},\citeyear{Tsybakov82b},\citeyear{Tsybakov83}] obtained
similar results to estimate the locally almost constant functions.

\bigskip

\paragraph{Minimax Estimation.} To guarantee good performance of the $ \rho $-LPA estimator in the minimax sense, we
assume that $
\rho' $ is bounded and Lipschitz. On the other hand, the
Huber function satisfies these assumptions, making it suitable for our problem. Moreover, it is commonly used in
practice (see for instance \cite{Petrus99} and \cite{Chang_Guo}). As for {\it linear
estimators} (kernel estimators, least square estimators, etc.), a good choice of the bandwidth $ h= \bar h_n(\beta,L) $
provides an optimal $\rho$-LPA estimator over the Hölder space $ \mathbb{H}_d({\beta},L,M) $. Finally,
$\bar h_n({\beta},L)=(L^2n)^{-\frac{1}{2\beta+d}}$ is chosen as the solution of the following bias/variance trade-off
\begin{equation}\label{equation_minimiser}
\big(nh^d\big)^{-1/2}+Lh^\beta\rightarrow\min_{h}. 
\end{equation}
In the model (\ref{model}), we show that the corresponding estimator
$\hat{f}^{\bar h_n({\beta},L)}({x_0})$ achieves the rate of convergence $ n^{-\beta/(2\beta+d)} $ (cf. Definition
\ref{def_minimax}) for $f({x_0})$ on
$\bH_d({\beta},L,M)$ (See Theorem \ref{minimax}).
We should point out that both the knowledge of $\beta$ and $L$ is required to the statistician in order to built the
optimal bandwidth $\bar h_n(\beta,L)$.

\bigskip

\paragraph{Adaptive Estimation.} In nonparametric statistics, an important problem is the adaptation compared to the
smoothness parameters $ \beta$ and $L $ that are unknown in practice. This
requests to develop a data-driven (adaptive) selection to choose the bandwidth. Then, the interesting
feature is the selection of estimators from a given family (cf. \cite{Barron_Birge_Massart99},
\cite{Lepski_Mammen_Spokoiny97}, \cite{Goldenshluger_Lepski08}). In this context, several approaches to the
selection from the family of linear estimators were recently proposed, see for instance \cite{Goldenshluger_Lepski08},
\cite{Goldenshluger_Lepski09}, \cite{Juditsky_Lepski_Tsybakov09} and the
references therein. However, those methods strongly rely on the linearity property. Robust estimators are generally
non-linear, there standard arguments (like the bias/variance
trade-off) cannot be applied straightforwardly. For instance, \cite{Brown_Cai_Harrison08} use
the asymptotic normality of the median to approximate the model (\ref{model}) by the
{\it wavelet sequence data} and they use {\it BlockJS wavelet thresholding} for adaptation over Besov spaces with the
integrated risk. Recently, 
\cite{Reiss_Rozenholc_Cuenod09} have considered the pointwise estimation for locally almost constant functions in
the homoscedastic regression with a heavy-tailed noise. That corresponds to $ \beta\leq1 $ for the Hölder functions in
the model (\ref{model}) (cf. also Definition \ref{def_holder}). They have considered the symmetric and continuous
density with $ g(0)>0 $. 
\medskip

In the context of adaptation, other new points are developed in this
paper:
\begin{itemize}
 \item adaptative pointwise estimation for any regularity $ \beta $ of isotropic
functions,
\item random design and heteroscedastic model,
\item unknown and heavy-tailed noise.
\end{itemize}
\medskip

For it, we construct an adaptive estimator (cf. Definition \ref{def_admissible}) using
general adaptation scheme due to \cite{Lepski90} ({\it Lepski's method}). This method is applied to choose the bandwidth
of the $ \rho $-LPA
estimator in the model (\ref{model}). 

We remind that $M$, the upper bound of $f$ and its partial derivatives, is involved in the construction of the
$ \rho $-LPA estimator (\ref{Intro_criterion}). Then, we assume that the parameter $ M $ is known and we do not study
the adaptation compared to
it. Contrary to the constants $ \beta,L $, one could estimate $ M $ to ``inject" it in the procedure without loss
of generality in the performance of our estimator (cf. \cite{Hardle_Tsybakov88}).

\bigskip

\paragraph{Exponential Inequality.}Lepski's procedure requires, in particular to
establish the exponential inequality for the deviations of $
\rho $-LPA estimator. As far as we know, these results seems to be new.

Denote by $ \bP_f $ the probability law of the observations $ Z^{(n)} $ satisfying (\ref{model}).
As we mentioned above, we need to establish the following inequality, for any $
\varepsilon>0\: $ and $ h\in(n^{-1/d},1) $:
\begin{equation}\label{expo_ineq}
 \mathbb{P}_f\left(\big|\hat{f}^h({x_0})-f({x_0})\big|\geq
\frac{\varepsilon}{\sqrt{nh^d}}\right)\leq
\cC\exp\left\{-\frac{\varepsilon^2}{
A+B\varepsilon/\sqrt{nh^d}}\right\},
\end{equation}
where $\cC,\:A,\:B$ are positive constants and $ A,\:B$ must be ``known". Details are given in Proposition \ref{lemma1}.
All results of this paper are based on
(\ref{expo_ineq}). 

The main difficulty in establishing (\ref{expo_ineq}) is that the explicit expression of $ \rho$-LPA estimator is not
typically available. Let us briefly discuss the main ingredients of M-estimation allowing to prove (\ref{expo_ineq}). If
the derivative of
contrast function $ \rho' $ is continuous, then solving the minimization problem (\ref{Intro_criterion}) can be viewed
as solving the following system of equations in $ t $ (first order condition):
\begin{equation}\label{Intro_partial derivative}
\forall p,\quad\tilde D_h^p(t):=\sum_{i:X_i\in V_{x_0}(h)}\frac{\partial}{\partial
t_p}\rho\big(Y_i-f_t(X_i)\big)\:K\left(\frac{X_i-{x_0}}{h}\right)=0,
\end{equation}
where $ t_p $ is the $ p^{th} $ component of the vector $ t $. Since $ \rho' $ is bounded, the partial
derivatives $ \tilde D_h^p(\cdot) $ can be viewed as an empirical process (i.e. a sum of independent and
bounded random
variables).

Denote $ \tilde D_h(\cdot) $ the vector of partial derivatives and $ D_h(\cdot)=\bE_{f_\theta}\tilde D_h(\cdot) $ where
$
\bE_{f_\theta}=\bE_{f_\theta}^n $ is the mathematical expectation with respect
to the probability law $ \bP_{f_\theta} $ of the ``parametric" observations $ (X_i,\cY_i)_{i=1,...,n} $.

\begin{figure*}[htb!]
\begin{minipage}[b]{0.99\linewidth}
\centerline{
\includegraphics[width=10cm]{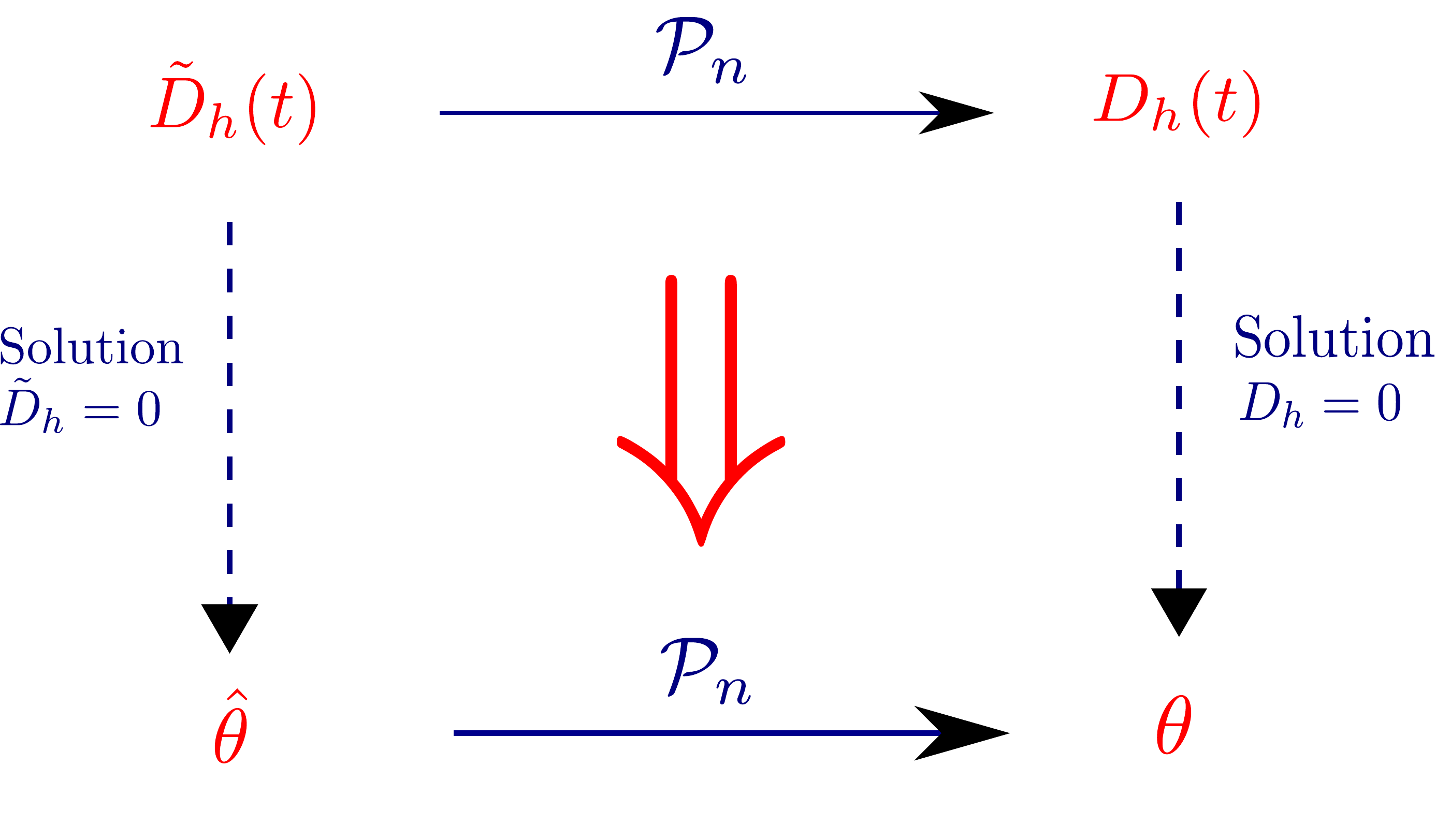}}
\end{minipage}
\caption{Illustration of the deviations' control. $ \cP_n $ represent the probability convergence.}
\label{fig_diagrammeMestimation}
\end{figure*}
Properties of the function $ D_h(\cdot) $ allow us to prove that $ \theta $ is the unique solution of $ D_h(\cdot)=0 $.
We
also notice
that $ |\hat{f}^h({x_0})-f({x_0})|\leq \|\hat\theta-\theta\|_1 $. The idea (presented in Figure
\ref{fig_diagrammeMestimation}) is to deduce the exponential inequality for
$\|\hat\theta-\theta\|_1$ from the exponential inequality for $ \sup_t\big|\tilde D_h(t)-D_h(t)\big| $. As we
mentioned above, we notice that $ \sup_t\big|\tilde D_h(t)-D_h(t)\big| $ can be viewed as the supremum of an
empirical process.

Now, classical arguments in probability tools can be used. To control $
\sup_t\big|\tilde D_h(t)-D_h(t)\big| $, we could
used standard tools developed by
\citeauthor{Talagrand96a} [\citeyear{Talagrand96a},
\citeyear{Talagrand96b}], \cite{Massart00} or
\cite{Bousquet02}.
But the obtained exponential inequalities (like
(\ref{expo_ineq})) contain unknown constants or require the
knowledge of an expectation's bound of $ \sup_t\big|\tilde D_h(t)-D_h(t)\big| $. To obtain this
bound, we can use the maximal inequalities developed by
\cite{VanderVaart_Wellner96}
(Chapter 2, Section 2.2) for {\it sub-gaussian processes}. 
But here again, there are universal constants (and unknown)
in the bound of the expectation.
\cite{Massart07} (Chapter 6) gives exponential inequalities for $ \sup_t\big|\tilde
D_h(t)-D_h(t)\big| $  without the expectation, but some constants are very big
in our case. In this paper, we choose to apply standard 
{\it chaining argument} and {\it Bernstein's inequality} (cf. (\ref{inequality_Bernstein}))
directly on $ \sup_t\big|\tilde D_h(t)-D_h(t)\big| $ (cf. Proof of Lemma \ref{large_deviation}). That allows us to have
constants smaller the ones cited in the papers above. 

\paragraph{Perspectives.}
\begin{itemize}
 \item We think that conditions on the noise densities could be reduced. We could consider the densities not necessary
monotonically on $ \bR_+ $, only the symmetric assumption seems necessary. 
\item A possible
perspective of this work is the study of estimating anisotropic functions. 
Indeed, the method developed by \cite{Kerkyacharian_Lepski_Picard01}, \cite{Klutchnikoff05} and
\citeauthor{Goldenshluger_Lepski08} [\citeyear{Goldenshluger_Lepski08}, \citeyear{Goldenshluger_Lepski09}] are based on
the linear properties and the machinery considered in those works can not adapt straightforwardly to nonlinear
estimators. 
\item Another perspective is to prove an oracle inequality for the family of $\rho$-LPA estimators indexed by
the bandwidth with the integrated risk. It could be interesting to introduce some criterion for choosing the optimal
contrast
function. 
\item Finally, we should also study
the heteroscedastic model (\ref{model}) with a degenerate design when the design density is vanishing or exploding.
\end{itemize}

\bigskip

This paper is organized as follows. We present exponential
inequalities in Section \ref{section_expo_ineq}, in order to
control deviations of $ \rho $-LPA
estimator. In Section \ref{sectionmaximalrisk}, we present
the results concerning minimax estimation and Section
\ref{sectionAdaptive} is
devoted to the adaptive estimation. An application of $
\rho $-LPA estimator with Huber function is proposed in
Section \ref{section:_huber}.
The proofs of the main results (exponential inequalities and
upper bounds) are given in Section
\ref{sectionUpper_bounds}, technical lemmas are postponed to the appendix.

\section{Exponential inequality for $ \rho $-LPA estimator}\label{section_expo_ineq}

\paragraph{Construction.}

To construct our estimator, we use the so-called {\it local
polynomial approach} (LPA) which consists in the following.
Let
$$
V_{x_0}(h)=\bigotimes_{j=1}^d\big[y_j-h/2,y_j+h/2\big]\cap[0,1]
^d ,
$$
be a neighborhood around $x_0$ of width $h\in(0,1)$. Fix $b>0$ (without loss of generality we will assume
that $b$ is an integer), let 
$$ \cS_b=\{  p=(p_1,\ldots,p_d)\in\bN^d:\;0\leq |p|
\leq b,\: |p|=p_1+\ldots+p_d\},
$$ 
and we denote $N_b$ the cardinal of $ \cS_b $.
Let $U(z), z\in\bR^{d}$ be the $N_b$-dimensional vector of
monomials of the following type (the sign $\top$ below
denotes the
transposition):
\begin{equation}\label{objet_polynome}
 U^\top (z)=\left(\prod_{j=1}^d
z_j^{p_j},\;\;  p\in\cS_b\right).
\end{equation}
For any $t^\top=\big(t_{p_1,\ldots,p_d}\in\bR\::\: 
p\in\cS_b\big)\in\bR^{N_b} $, we define the local
polynomial in a neighborhood of $x_0$ as for any $ x\in[0,1]^d$
\begin{eqnarray}\label{polynome}
f_t(x)&=&t^\top
U\left(\frac{x-{x_0}}{h}\right)\bI_{V_{x_0}(h)}(x)=\sum_{p\in\cS_b}
t_p\:\left(\frac { x-{x_0}}{h}\right)^p\bI_{V_{x_0}(h)}(x),
 \end{eqnarray}
where $ z^p=z_1^{p_1}\cdots z_d^{p_d} $ for $
z=(z_1,\dots,z_d) $ and $\bI$ denotes the indicator
function.
For any $ M>0 $, introduce the following subset
of $\bR^{N_b}$
\begin{eqnarray}
\label{ensemble} \Theta\big(M\big)=\left\{t\in\bR^{N_b}:\;\;
\|t\|_1\leq M\right\},
\end{eqnarray}
where $\|\cdot\|_1$ is  $\ell_1$-norm on $\bR^{N_b}$.
We notice that for any $
t\in\Theta(M)
$, $ \|f_t\|_\infty\leq M $  where $\|\cdot\|_\infty$ is  the sup-norm on $[0,1]^d$.

\medskip
The function $ \rho $ is called {\it contrast function} if
it has the following properties.
\begin{assumptions}\label{properties_rho}
\hspace{1cm} 
\begin{enumerate}
 \item $ \rho : \bR\rightarrow\bR_+ $ is symmetric, convex and $ \rho(0)=0 $,
 \item the derivative $ \rho' $ is $ 1 $-Lipschitz on $ \bR
$ and bounded: $ \dot{\rho}_\infty=\|\rho'\|_\infty<\infty $,
\item the second derivative $ \rho'' $ is defined almost
everywhere and there exist $L_{\rho}>0$ and $\alpha>0$ such that
$$
\sup_{i=\overline{1,n}}\int_\bR\big|\rho''(z+u)-\rho''(z+v)\big|\:g_i(z)\:dz\leq
L_{\rho}\:|u-v|^\alpha,\quad \forall u,v\in\bR.
$$
where $\overline{1,n}=1,\ldots,n$.
\end{enumerate}
\end{assumptions}
A well-known example of a contrast function $ \rho $ satisfying Assumption \ref{properties_rho} above is the Huber
function (cf.
\cite{Huber64}) presented in Section \ref{section:_huber}.

\medskip

Let $ K $ be a kernel function, i.e. a positive function with a compact support included in $ [-1/2,1/2]^d $ such that $
K_\infty:=\|K\|_\infty<\infty $ and $ \int_\bR K(z)dz=1 $.
We will construct the $ \rho $-LPA estimator for $f({x_0})$ using {\it local $ \rho $-criterion} which is defined as
follows:
\begin{eqnarray}
\label{rho_criterion} 
\tilde\pi_h(t)=\tilde\pi_h\big(t,Z^{(n)}\big)=\frac{1}{nh^d}
\sum_{i=1}^{n}\rho\big(Y_i-f_t(X_i)\big)\:K\left(\frac{X_i-{x_0}
}{h}\right).
\end{eqnarray}
Let $ \hat\theta(h) $ be the solution of the following
minimization problem:
\begin{equation}
 \label{minimization}
   \hat\theta(h) = \arg\min_{t\in\Theta(M)}\tilde\pi_h(t).
 \end{equation}
The {\it $ \rho $-LPA estimator} $ \hat{f}^{h}({x_0}) $ of $ f({x_0}) $ is defined as $\hat{f}^{h}({x_0}) =
\hat\theta_{0,...,0}
(h)$. We notice that this local approach can be considered as the estimation for successive derivatives of the function
$ f $. However in the present paper, we focus on the estimation of $ f({x_0}) $.

%


\paragraph{Exponential inequality.} Later on, we will only consider values of $h\in[n^{-1/d},1]$. Put
$\theta=\theta(f,x_0,h)=\big\{\theta_{  p}:\; 
p\in\cS_b\big\}$, where
$\theta_0=\theta_{0,...,0}=f({x_0})$  and
\begin{equation}\label{coefficient_taylor}
\forall p\in\cS_b\::\:p\neq0,\quad\theta_{  p}=\displaystyle \frac{\partial^{| 
p|}f({x_0})}{\partial y_1^{p_1}\cdots\partial
y_d^{p_d}}\frac{h^{| 
p|}}{p_1!...p_d!}.
\end{equation}
Here, we do not assume the existence of partial derivatives of $f$. To define $ \theta $ properly, the following
agreement will be used in the sequel:
if the function $f$ and the vector $ {p}$\; are such that $\partial^{| 
p|}f$ does not exist, we put $\theta_{p}=0$.



Set $\cB\big(\theta,z\big)=\left\{
t\in\Theta(M):\:\|t-\theta\|_2\leq z\right\}$ the Euclidean ball with radius $z$ and center $\theta$ and define the
event for any $ h,z>0 $
\begin{eqnarray}\label{voisinageSet}
G_z^h&=&\big\{\hat\theta(h)\in\cB\big(\theta,z\big)\big\},\quad
\end{eqnarray}
where $ \hat\theta(h) $ is given by (\ref{minimization}). Let
\begin{equation}\label{def_constant_sigma}
\Sigma=2+2\sum_{l=1}^\infty d^210^{2l-1}
\:\exp\left\{-\frac{18\:10^l}{\pi^4\:l^4}\frac{(8K_\infty)^{-1}}{K_\infty+1/3}
\right\},
\end{equation}
be some finite constants
and let the constant $ \lambda $ be the smallest eigenvalue of matrix
$$ 
\int_{[-1/2,1/2]^d}U(x)
\:U^\top(x)\:K(x)\:dx.
$$
\cite{Tsybakov08} (Lemma 1.4) showed that $ \lambda $ is positive, on the hand the last matrix is strictly positive
definite.
Finally, put
\begin{equation}\label{def_constant_c}
 c\big(\rho,(g_i)_i\big)=\inf_{i=1,\ldots,n}\int_\bR \rho''(z)\:g_i(z)\:dz,
\end{equation}
and define the set of sequences of symmetric densities which are monotonically on $ \bR_+ $
\begin{equation}\label{def_set_densities}
 \cG_{\rho}^{(c)}=\left\{(g_i)_i\::\:c\big(\rho,(g_i)_i\big)\geq c\right\},\quad c>0.
\end{equation}
Denote for all $ a,b\in\bR,\: a\vee b=\max(a,b) $. The next proposition is the  milestone for all results
proved in the
paper.

\begin{proposition}\label{lemma1} Let $ \rho $ be a contrast function and let $ c>0 $. Then, for any $ n\in\bN^* $, $
(g_i)_i\in\cG_{\rho}^{(c)}
$,
$ h>n^{-1/d} $, $ M>0 $, $ x_0\in[0,1]^d $ and any $ f $ such that $ \|\theta(f,x_0,h)\|_1\leq M $, 
we have for any $ \varepsilon\geq\frac{4
N_b}{{c\lambda}}\big(1\vee b_h\:\sqrt{nh^d}\big)$
\begin{eqnarray}\label{eq_lemma1}
&&\mathbb{P}_f\left(\sqrt{nh^d}\big|\hat{f}
^h({x_0})-f({x_0})\big|\geq\varepsilon,\:G_\delta^h\right)\nonumber\\
&&\quad\leq
N_b\varSigma\exp\left\{-\frac{\left(\frac{
c\lambda}{{2N_b}}\:\varepsilon
-\big(1\vee b_h\:\sqrt{nh^d}\big)\right)^2}{8
K_\infty^2(1\vee\dot{\rho}_\infty^2)+\frac{4\: K_\infty}{3{N_b}}\:(1\vee{\dot{\rho}_\infty})\frac{
c\lambda\:\varepsilon}{
\sqrt{nh^d}} }\right\}.
\end{eqnarray}
\end{proposition}
The proof of this proposition is given in Section \ref{sectionUpper_bounds}.
\begin{remark}
 The control of the deviations of $ \hat f^h $ is realized under the event $ G_\delta^h $ that the estimator $
\hat\theta(h) $ is contained in a ball centered at $ \theta $ whereas its radius does not depend on $ n $, else
 it could change the rate of convergence. In Section \ref{sectionUpper_bounds} we give an
exponential inequality to control the probability of the complementary of $ G_\delta^h $ (cf.
Lemma \ref{Control_complement}).
\end{remark}

\begin{remark}
In the minimax case, the knowledge of constants in (\ref{eq_lemma1}) is not required. However for adaptation, the
constant $ c $ is involved in the construction of adaptive estimator. This restricted the consideration of the noise
densities which satisfy (\ref{def_set_densities}). We notice that this problem is simplified to the calibration of an
alone constant with a dataset.
\end{remark}

\section{Minimax Results on $\mathbb{H}_d({\beta},L,M)$}\label{sectionmaximalrisk}

In this section, we present several results concerning
maximal and minimax risks on $\mathbb{H}_d({\beta},L,M)$. We propose the estimator which bound the
maximal risk on this class of
functions without restriction imposed on these parameters.

\paragraph{Preliminaries.}

\begin{definition}
\label{def_holder} Fix $\beta>0$, $L>0$, $M>0$ and  let
$\lfloor\beta\rfloor$ be the largest integer  strictly smaller
than
$\beta$. The \emph{isotropic Hölder class}
$\bH_d({\beta},L,M)$ is
the set of functions $f:[0,1]^d\rightarrow\bR$ admitting on
$[0,1]^d$
all partial derivatives of order $\lfloor\beta\rfloor$ and
such that for any
$ x,y\in [0,1]^d$

\begin{eqnarray*}
\left|\frac{\partial^{|p|}f(x)}{\partial
x_1^{p_1}\cdots\partial x_d^{p_d}}-\frac{\partial^{| 
p|}f({y})}{\partial y_1^{p_1}\cdots\partial
y_d^{p_d}}\right|&\leq&L
\big[\|x-{y}\|_1\big]^{\beta-\lfloor\beta\rfloor},
\quad\forall\;
\big| {p}\big|=\lfloor\beta\rfloor,\\*[2mm]
\sum_{p\in\cS_{\lfloor\beta\rfloor}}\sup_{x\in[0,1]^d}
\left|\frac{\partial^{| 
p|}f(x)}{\partial
x_1^{p_1}\cdots\partial x_d^{p_d}}\right|&\leq& M.
\end{eqnarray*}

where $x_j$ and $y_j$ are the $j^{th}$ components of $x$ and
$y$.

\end{definition}

Let $\bE_f=\bE^{n}_f$ be the mathematical expectation with
respect
to the probability law $ \bP_f $ of the observation
$Z^{(n)}$ satisfying
(\ref{model}). Firstly, we define the maximal risk on
$\mathbb{H}_d({\beta},L,M)$ corresponding to the estimation
of the
function $f$ at a given point $x_0\in[0,1]^{d}$.

Let $\tilde{f}_n$ be an arbitrary estimator built from the
observation $Z^{(n)}$. For any $r$ let
\begin{equation}\label{def_maximal_risk}
 R_{n,r}\big[\tilde{f}_n,\mathbb{H}_d({\beta},L,M)\big]
=\sup_
{f\in\mathbb{H}_d({\beta},L,M)}\mathbb{E}_f\big|\tilde{f}
_n({x_0})-f({x_0})\big|^r.
\end{equation}
This quantity is
called
{\it maximal risk} of the estimator $\tilde{f}_n$ on
$\mathbb{H}_d({\beta},L,M)$ and the {\it minimax risk} on
$\mathbb{H}_d({\beta},L,M)$ is defined as
\begin{equation}\label{def_minimax_risk}
 R_{n,r}\big[\mathbb{H}_d({\beta},L,M)\big]=\inf_{\tilde{f}}R_{n,r}\big[\tilde{f},\mathbb{H}_d({\beta},L,M)\big],
\end{equation}
where the infimum is taken over the set of  all 
estimators.
\begin{definition}
\label{def_minimax} The normalizing sequence $\psi_n$ is called
minimax rate of convergence and the estimator $\hat{f}$ is
called minimax (asymptotically minimax) if
\begin{eqnarray}
\label{def_lower_bounds}\liminf_{n\to\infty}\psi^{-r}_n\:R_{n,r}\big[\hat{f},\mathbb{H}_d({\beta},L,M)\big]&>&0,\\
\label{def_upper_bounds}\limsup_{n\to\infty}\psi^{-r}_n\:R_{n,r}\big[\hat{f},\mathbb{H}_d({\beta},L,M)\big]&<&\infty.
\end{eqnarray}
\end{definition}

\paragraph{Upper bound for maximal risk.}
Let the minimizer of the bias/variance trade-off
(\ref{equation_minimiser}) be given by
\begin{equation}
\label{hbar} \bar h=(L^2n)^{-\frac{1}{2\beta+d}}.
\end{equation}
 
The next theorem shows how to construct the estimator based
on locally parametric approach which achieves the following rate of convergence in the model (\ref{model})
 \begin{equation}
 \label{lower_sequence}
 \varphi_n({\beta})=n^{-\frac{{\beta}}{2{\beta}+d}}.
 \end{equation}
Let $\hat{f}^{\bar
h}({x_0})=\hat{\theta}_{0,\ldots,0}\big(\bar h\big)$ be given by
(\ref{ensemble}), (\ref{rho_criterion}) and 
(\ref{minimization}) with
$h=\bar{h}$ and $b=\lfloor\beta\rfloor$.

\medskip

\begin{theorem}\label{minimax}
Let $ \beta>0 $, $ L>0 $, $ M>0 $, $ x_0\in[0,1]^d $, $ c>0 $ and let $\rho$ be a fixed contrast function. Then
for any $
(g_i)_i\in\cG_{\rho}^{(c)} $
\begin{eqnarray*}
\limsup_{n\to\infty}\:\varphi_n^{-r}(\beta)\:R_{n,r}\Big[
\hat{f}^{\bar
h}({x_0}),\mathbb{H}_d({\beta},L,M)\Big]<\infty,\quad \forall
r\geq1.
\end{eqnarray*}
\end{theorem}
This theorem will be deduced from Proposition
\ref{lemma1} and the proof is given in Section \ref{section_proof_maximalrisk}.
\begin{remark}\label{remark_minimax}
\cite{Tsybakov82a} showed lower bounds (\ref{def_lower_bounds}) for rate
$n^{-\frac{\beta}{2\beta+d}}$ with the following assumption on Kullback distance on
the noise density $ g $, i.e. it
exists $ v_0>0 $ such that
$$
\int g(u)\:\ln\frac{g(u)}{g(u+v)}\:du\leq
o\big(v^2\big),\quad\forall v\::\:|v|\leq v_0.
$$
We notice that Gaussian and Cauchy densities verify this
assumption (cf. also \cite{Tsybakov08} Chapter 2). In this case, we conclude that $ \hat f^{\bar h} $
is minimax and $
\varphi_n(\beta) $ is the minimax rate on
$\mathbb{H}_d({\beta},L,M)$.
\end{remark}

\section{Bandwidth Selection of $\rho$-LPA Estimator}\label{sectionAdaptive}

This section is devoted to the adaptive estimation over the
collection of classes
$\Big\{\mathbb{H}_d(\beta,L,M)\Big\}_{\beta,L}$. Here we suppose
$ M $ known, as we mentioned in the introduction, the parameter $ M $ could be estimated and used with a ``Plug-in"
method (cf. \cite{Hardle_Tsybakov92}). We will not impose any restriction on the possible value of $L$, but we
will assume that $\beta\in (0,b]$, where $b$ as previously, is an arbitrary chosen integer.

\medskip

We start by remarking that there is not optimally adaptive estimator. Well-known disadvantage of maximal approach is the
dependence of the estimator on the smoothness parameters
describing the functional class on which the maximal risk
is determined (cf. (\ref{def_maximal_risk})). In
particular, 
$\bar h_n({\beta},L)$, optimally chosen in view of
(\ref{equation_minimiser}), depends explicitly on ${\beta}$
and $L$. To overcome this drawback, a maximal adaptive
approach has been proposed by \cite{Lepski90} for
pointwise
estimation. The first question arising in the adaptation
(reduced to the problem at hand) can be formulated as
follows.

\smallskip

{\it  Does there exist an estimator which would be minimax
on
$\bH({\beta},L,M)$ simultaneously for all values of
${\beta}$ and $L$ belonging to some given set
$\mB\subseteq\big[\bR_{+}\setminus{0}\big]\times\big[\bR_{+}
\setminus{0}\big]$
?}

\smallskip

For integrated risks, the answer is positive (cf.
\cite{Lepski91},
\cite{Donoho_Johnstone_Kerkyacharian_Picard95},
\cite{Lepski_Spokoiny97},
\cite{Goldenshluger_Nemirovski97} and \cite{Judistsky97}).
For the estimation of the function  at a given point, it is
typical that the price to pay
is not null (cf. \cite{Lepski90}, 
\cite{Brow_Low96},
\cite{Lepski_Spokoiny97}, \cite{Tsybakov98},
\cite{Klutchnikoff05}, \cite{Reiss_Rozenholc_Cuenod09}, \cite{Chichignoud10}).
Mostly, the price to pay is a power of $ (b-\beta)\ln n $ for
pointwise estimation.

Let $ \Psi=\left\{\psi_n(\beta)\right\}_{\beta\in
(0,b]} $ be a given family of normalizations.
\begin{definition}\label{def_admissible}
The family $\Psi$ is called admissible  if there exists an
estimator
$\hat{f}_n$ such that for some $L,M>0$
\begin{equation}
\label{eq:def-admissibility} \limsup_{n\rightarrow\infty}
\psi^{-r}_n(\beta)\;R_{n,r}\big(\hat{f}_n,\mathbb{H}_d(\beta
,L,M)\big)<\infty,\;\;\forall
\beta\in (0,b].
\end{equation}
The estimator $\hat{f}_n$ satisfying
(\ref{eq:def-admissibility}) is
called $\Psi$-\textsf{attainable}. The estimator $ \hat{f}_n
$ is
called $\Psi$-\textsf{adaptive} if
(\ref{eq:def-admissibility})
holds  for \textsf{any} $ L>0 $.
\end{definition}

\cite{Lepski90} showed that the family of
rates $\left\{\varphi_n(\beta)\right\}_{\beta\in (0,b]}$, defined in (\ref{lower_sequence}), is not
admissible in the white noise model. With other tools, \cite{Brow_Low96} extend this result for density estimation and
nonparametric
Gaussian regression. It means that there
is no-estimator which would be minimax simultaneously
for several
values of parameter $\beta$, for pointwise estimation, even
if $L$ is supposed to be fixed. This result does not require
any
restriction on $\beta$ as well.

\medskip
Now, we need to find another family of normalizations for
maximal
risk which would be attainable  and, moreover, optimal in
view of
some criterion of optimality.
Let $\Phi$ be the following family of
normalizations, for any $\beta\in
(0,b]$
\begin{equation}\label{normalizationAdaptive}
 \phi_{n}(\beta)=\left(\frac{\varrho_{n}(\beta)}{n}\right)^{
\frac{\beta}{2\beta+d}},\quad
\varrho_{n}(\beta)=1+\frac{2(b-\beta)}{
(2\beta+d)(2b+d)}\ln n.
\end{equation}
We notice that $\phi_n(b)=\varphi_n(b)$ and for $n$ large enough
$\varrho_{n}(\beta)\sim{(b-\beta)\ln
n}$ for any $\beta\neq b$.
It is possible to show that this family $ \Phi $ is adaptive
optimal using the most recent criterion developed by
\cite{Klutchnikoff05} used for the white noise
model and used by \cite{Chichignoud10} for the multiplicative uniform regression. On the other hand, the so-called {\it
price to pay for adaptation} $\varrho_{n}(\beta)$ could be considered as optimal. 

\paragraph{Construction of $\Phi$-adaptive estimator.}We
begin by stating that the construction of our estimation
procedure is decomposed in several steps. First, we
determine the
family of $\rho$-LPA estimators. Next, based on Lepski's method, we propose a
data-driven selection from this family.

\medskip

Let $\rho$ be a fixed contrast function. In the model (\ref{model}), we recall that the sequence
of densities $ (g_i)_i $ is ``unknown" for the statistician. We take $
\hat{f}^{h} $ the estimator given by (\ref{ensemble}),
(\ref{rho_criterion}) and (\ref{minimization}), so the
family of $ \rho $-LPA estimators $\hat{\cF}$ is defined now
as follows. Put
\begin{equation}\label{def_hmin_hmax}
h_{\min}=(\ln n)^{2/d} n^{-1/d},\quad h_{\max}=n^{-\frac{1}{2b+d}},
\end{equation}
and
$$
h_k=2^{-k}h_{\max},\;\; k=\overline{0, \mathbf{k}_n}:=0,\ldots,\mathbf{k}_n,
$$
where $\mathbf{k}_n$ is the largest integer such that
$\displaystyle
h_{\mathbf{k}_n}\geq
h_{\min}$. Set
\begin{eqnarray}\label{eq:family-estimators}
\hat{\cF}=\left\{\hat{f}^{(k)}({x_0})=\hat{\theta}_{0,\ldots,0}
(h_k),\;\;k=\overline{0, \mathbf{k}_n}\right\}.
\end{eqnarray}
We put $\hat{f}^*({x_0})=\hat{f}^{(\hat{k})}({x_0})$, where
$\hat{f}^{(\hat{k})}({x_0})$ is
selected from $\hat{\cF}$ in accordance with the rule:
\begin{eqnarray}\label{indexe adaptive}
\hat
k=\inf\left\{k=\overline{0,\mathbf{k}_n}
:\;\;\big|\estimatorb{f}^{(k)}({x_0})-\hat{f}^{(l)}({x_0})\big|\leq
CS_n\big(l\big),\;\;
l=\overline{k+1,\mathbf{k}_n}\right\}.
\end{eqnarray}
Here we have used the following notations. Let $ c>0 $ be fixed and
\begin{eqnarray}\label{variance}
C&=&{\frac{4N_b}{c\lambda}}\left(1+2\:K_\infty\:(1\vee \dot{\rho}_\infty)\sqrt{rd}\right),\\
S_n(l)&=&\left[\frac{1+l\ln
2}{n\big(h_l\big)^d}\right]^{1/2},\quad
l=\overline{0, \mathbf{k}_n},\nonumber
\end{eqnarray}
where $ r\geq1 $ is the power of the risk and $ c $
is defined in (\ref{def_constant_c}), $ {\dot{\rho}_\infty} $ and $ K_\infty $ are respectively bounds of $
\rho'(\cdot) $ and $ K(\cdot) $, and the positive constant $ \lambda$ is the
smallest eigenvalue of the matrix $ \int_{[-0.5,0.5]^d}U(x) \:U^\top(x)
\:K(x)\:dx $. We will see that this matrix is strictly positive definite (cf. Lemma \ref{criterion_inversible}).

\medskip

\paragraph{Main Result.} The next theorem is the main result of this paper. It
allows us to guarantee a good performance of our adaptive $
\rho $-LPA estimator $ \hat f^* $.

\begin{theorem}
\label{lepski} Let $b>0, M>0$ and $ \rho $ be a fixed
contrast function. Then, for any $ (g_i)_i\in\cG_{\rho}^{(c)}
$, $\beta\in (0,b]$, $ L>0 $ and $r\geq1$
$$
\limsup_{n\to\infty}\phi_{n}^{-r}(\beta)\;R_{n,r}\Big[\hat{f
} ^*({x_0}),\mathbb{H}_d({\beta},L,M)\Big]<\infty.
$$
\end{theorem}
The proof (given in Section \ref{section_proof_lepski}) is based on the scheme due to \cite{Lepski_Mammen_Spokoiny97}.
\begin{remark}
The assertion of the theorem means that the proposed
estimator
$\hat{f}^*({x_0})$ is $\Phi$-adaptive in the model (\ref{model}) (cf. Definition
\ref{def_admissible}). It implies in particular that the
family of normalizations $\Phi$ is admissible. 
\end{remark}


\begin{remark}
In the present paper, we do not give the explicit expression of the constant in the upper bound
of the risk with the proof given in this paper. But it is possible to solve this problem. In the proof
of Lepski's method, we notice that the upper bound polynomially depends on the parameter $ C $ and it is important to
minimize this
constant. We see that this constant depends on the contrast function $ \rho $ and it is easy to see that minimizing $
C=C(\rho) $
can be viewed as minimizing the following Huber variance (cf. \cite{Huber_Ronchetti09} Page 74)
$$
\sigma^2_\rho=\frac{\int\big(\rho'\big)^2\:dg}{
\left(\int\rho''\:dg\right)^2}\rightarrow\min_{\rho},
$$
where $ g $ is the noise density in the homeoscedastic model.
\end{remark}

\begin{remark}
The limitation concerning the
consideration of isotropic
classes of functions is due to the use of Lepski's
procedure. It seems that to be able to treat the adaptation
over the scale of anisotropic
classes (i.e. $d$-dimensional functions with different regularities $ \beta $ for each variable). Another scheme should
be applied as in
\cite{Lepski_Levit99}, \cite{Kerkyacharian_Lepski_Picard01},
\cite{Klutchnikoff05} and
\cite{Goldenshluger_Lepski08}. As we have mentioned above,
these latter procedures cannot be used with $ \rho $-LPA
estimators, and for the model (\ref{model}) this problem is
still open.
\end{remark}

\section{Application: Huber function}\label{section:_huber}
Consider the model (\ref{model}), with following additional
assumptions.
\begin{equation}\label{eq:g_hetero}
 g_i(\cdot)=g(\cdot/\sigma_i)/\sigma_i,\quad i=1,...,n,
\end{equation}
where the density $ g $ is symmetric and monotonically on $ \bR_+ $. $ (\sigma_i)_i $ is
a sequence of real values such that for any $
 i,\:0<\sigma_{\min}\leq\sigma_i<\infty $ where $
\sigma_{\min} $ is known. The model (\ref{model})
with (\ref{eq:g_hetero}) can be written as 
\begin{equation}\label{model_exemple}
 Y_i=f(X_i)+\sigma_i\:\xi_i,\quad i=1,...,n,
\end{equation}
where $ (\xi_i) $ are i.i.d. with the density $g$.
%

Let
\begin{equation}\label{huber_function}
  \rho_\gamma(z) = \gamma(z-0.5\:\gamma)\:\bI_{|z|>\gamma}+0.5\:z^2\:
\bI_{|z|\leq\gamma},\:z\in\bR,\:\gamma\geq 0.
\end{equation}
the Huber function (\cite{Huber64}).
We construct the $ \rho_\gamma $-LPA estimator from (\ref{ensemble}),
(\ref{rho_criterion}) and (\ref{minimization}). The function $ \rho_\gamma $ is a {contrast function} verifying
Assumption
\ref{properties_rho}. Recall that the constant $ c=c\big(\rho_\gamma,(g_i)_i\big) $ defined in
(\ref{def_constant_c}) must be positive. We notice
that the second derivative can be written as $
\rho_\gamma''(\cdot)=\bI_{[-\gamma,\gamma]}(\cdot) $ and that
$$
c\big(\rho_\gamma,(g_i)_i\big)\geq c_\gamma:=2\int_{0}^{
\gamma\sigma_{\min}} g(z)dz.
$$
We formulate the following assertion: for any $
\sigma_{\min}>0 $ and any $ g $ a symmetric density and monotonically on $ \bR_+ $, there exists a
constant $ \gamma_0>0 $ such that for any $ \gamma\geq\gamma_0,\: c_{\gamma}>0 $.

We propose the adaptive $ \rho_{\gamma_0} $-LPA
estimator $ \hat f^*_{\gamma_0}({x_0}) $ selected with the data-driven selection proposed in
Section \ref{sectionAdaptive} with the constant
$$
C={\frac{2N_b}{\lambda\int_{0}^{
\gamma_0\sigma_{\min}} g(z)dz}}\left(1+2\:K_\infty\:(1\vee \gamma_0)\sqrt{rd}\right).
$$
The next result is a direct
consequence of Theorem \ref{lepski}.

\begin{corollary}\label{corollary_lepski}
 Let $b>0$, $ M>0$ be some fixed constants and consider the model
(\ref{model_exemple}). Then, for any $ (g_i)_i\in\cG_{\rho_{\gamma_0}}^{(c_{\gamma_0})}
$, $\beta\in (0,b]$, $ L>0
$ and $r\geq1$
$$
\limsup_{n\to\infty}\phi_{n}^{-r}(\beta)R_{n,r}\Big[\hat{f}
^*_{\gamma_0}({x_0}),\mathbb{H}_d({\beta},L,M)\Big]<\infty.
$$
\end{corollary}
\begin{remark}
 We notice that the threshold parameter $ C $ explicitly depends on the minoration $ \sigma_{\min} $ of the noises
variances $ (\sigma_i)_i $. Contrary to linear estimators (C polynomially depends on $(\sigma_i)_i$), we can see that
the influence of $ (\sigma_i)_i $ is very limited for $ \rho $-LPA estimators.
\end{remark}

\begin{remark}\label{Remark_estimator_semirobust}
Corollary \ref{corollary_lepski}
only guarantees that asymptotically for any $ \gamma\geq\gamma_0 $, $ \rho_\gamma $-LPA estimators have the same
performance.
In the future, an important question  to adress is: how one can choose the parameter $ \gamma $? In
theory, there is yet no criterion for choosing an optimal $ \gamma $, but we can make the following remarks. If $
\gamma=\infty
$, then the $ {\rho}_\infty $-LPA estimator is the least square estimator (sensitive to
extreme values of the noise) and if $ \gamma=0 $ then the $ \rho_0 $-LPA estimator becomes the median estimator (robust
estimator). 
It is well-known that least squares estimator and median estimator respectively suffer from undersmoothing and
oversmoothing. This phenomenon is highlighted by
\cite{Reiss_Rozenholc_Cuenod09}. We believe that a better choice of parameter $ \gamma $ should give a
``semi-robust" estimator.
Locally this could reduce the above mentioned issue. In
practice, it will be interesting to
select
the parameter $ \gamma $ as a measurable function of observations which adapts to extreme values of the noise. This
problem is related to the
estimation of the noise variance and to the minimization of the Huber variance (cf. \cite{Huber_Ronchetti09} Page 74).
\end{remark}

\section{Proofs of Main Results: Exponential inequalies and
Upper Bounds}\label{sectionUpper_bounds}

\subsection{Proof of Proposition \ref{lemma1}}

%

\paragraph{Notations.}
Recall that $ \cS_b=\big\{  q\in\bN^d\::\:|  q|\leq
b,\:  |  q|=q_1+\ldots+q_d \big\}$ and $ N_b $ its cardinal. We consider the {\it partial derivative} of the local $
\rho
$-criterion
\begin{eqnarray}\label{Derivate_rho_criterion}
\tilde D_h(\cdot)=\left(\frac{\partial}{\partial
t_{p}}\:\tilde\pi_h(\cdot)\right)_{p\in\cS_b}^\top,
\end{eqnarray}
where $\tilde\pi_h(\cdot)$ is the local $ \rho $-criterion defined in (\ref{rho_criterion}).
Let also
\begin{equation}\label{Deterministic_criterion}
 \cE_h(\cdot)=\bE_f\left[\tilde D_h(\cdot)\right] \quad\text{and}\quad D_h(\cdot)=\bE_{f_\theta}\big[\tilde
D_h(\cdot)\big],
\end{equation}
where $ f_\theta $ is the Taylor polynomial defined in
(\ref{coefficient_taylor}), $ \bE_{f_\theta}=\bE_{f_\theta}^n $ be the mathematical expectation with respect to
the probability law $ \bP_{f_\theta} $ of the ``parametric" observations $ (X_i,\cY_i)_{i=1,...,n} $ (cf.
(\ref{model_parametric})) and $ \bE_{f}=\bE_{f}^n $ be the mathematical expectation with respect
to the probability law $ \bP_{f} $ of the observation $Z^{(n)}$.

We call the {\it Jacobian matrix} $ J_D $ of $ D_h $ such that
\begin{eqnarray}\label{Jacobian_matrix}
\big(J_D(\cdot)\big)_{ 
p, 
q\in\cS_b}:=\left(\frac{\partial }{\partial t_{q}}D_h^{ 
p}(\cdot)\right)_{  p, 
q\in\cS_b}=\left(\bE_{f_\theta}\frac{\partial^2}{\partial t_{p}\partial t_{q}}\:\tilde \pi_h(\cdot)\right)_{ 
p, 
q\in\cS_b},
\end{eqnarray}
where $ D_h^{p}(\cdot) $ is the $ p^{th} $ component of $ D_h(\cdot) $.

\paragraph{Auxiliary lemmas.}
We give the following lemma concerning the {\it
deterministic criterion} $ D_h $ defined in
(\ref{Deterministic_criterion}).
Denote $ \|\cdot\|_2 $ the $ \ell_2 $-norm on $\bR^{N_b}$.
\begin{lemma}\label{criterion_inversible}
Let $ \rho $ be a contrast function, for any $ (g_i)_i\in\cG_{\rho}^{(c)} $
we have the following assertions:
\begin{enumerate}
 \item the matrix $ J_D(\theta) $ is strictly positive definite and $ \theta $ is the unique
solution of the equation $ D_h(\cdot) = (0,...,0) $ on $ \Theta(M) $,
\item there exists $ \delta>0 $ which only depends on the contrast function $\rho$ such that for any $
\tilde\theta\in\cB(\theta,\delta) $, we have
$$
\|\tilde\theta-\theta\|_2\leq{\frac{2}{c\lambda}}
\inf_{h>0}\left\|D_h\big(\tilde\theta\big)-D_h(\theta)\right\|_2.
$$
\end{enumerate}

\end{lemma}
Recall that $ b_h=\sup_{x\in V_{x_0}(h)}\big|f_\theta(x)-f(x)\big| $
corresponds to the approximation error (bias) and denote $ \cE_h^{p}(\cdot) $ the $ p^{th} $ component of $ \cE_h(\cdot)
$.
Let us give a lemma which
allows us to control the {\it bias term}.

\begin{lemma}\label{control_bias}
For any contrast function $ \rho $, $ h>n^{-1/d} $ and any $ f $ such that $ \|\theta\|_1\leq M $,
we have
$$
\max_{p
\in\cS_b}\sup_{t\in\Theta(M)}\big|\cE_h^{  p}(t)-D_h^{  p}(t)\big|\leq b_h.
$$

\end{lemma}

The next result allows us to control deviations of
{\it partial derivatives of $ \rho $-criterion} $ \tilde D_h
$ defined in
(\ref{Derivate_rho_criterion}).

\begin{lemma}\label{large_deviation}
 For any contrast function $ \rho $, any $ f $ such that
$ \|\theta\|_1\leq M $ and any $ h>n^{-1/d} $, we have for any $ z\geq 2\big(1\vee b_h\:\sqrt{nh^d}\big) $
\begin{eqnarray*}
&&\max_{p\in\cS_b}\bP_f\left(\sqrt{nh^d}\sup_{t\in\Theta(M)}\left|\tilde
D_h^{  p}(t)-\cE_h^{  p}(t)\right|\geq z\right)\\
&&\quad\leq\varSigma\exp\left\{-\frac{\left(z-
b_h\:\sqrt{nh^d}\right)^2}{8
K_\infty^2\:(1\vee\dot{\rho}_\infty^2)+\frac{4K_\infty}{3\sqrt{nh^d}}\:(1\vee{\dot{\rho}_\infty})z}\right\}.
\end{eqnarray*}
\end{lemma}
As we mentioned above, the partial derivatives $ \big(\tilde D_h^p(\cdot)\big)_p $ can be considered as empirical
processes. Thus the proof (given in Appendix) is based on a chaining
argument and Bernstein's inequality (cf. (\ref{inequality_Bernstein})).
In particular, it is required that the derivative $ \rho'
$ of the contrast function is bounded and Lipschitz.

Denote by $ \bar G_\delta^h $ the complementary of $
G_\delta^h $ (defined in (\ref{voisinageSet})) where the radius $ \delta $ is defined in
Lemma \ref{criterion_inversible} and let $ \varkappa_\delta=\inf_{
t\in\Theta(M)\backslash\cB(\theta,\delta)}\|D_h(t)\|_2
/2 $ be a positive constant. The next lemma allows us to control
the probability of the event that ``the $ \rho $-LPA estimator
does not belong to the
ball  centered on $ \theta $ with radius $ \delta $".

\begin{lemma}\label{Control_complement}
 For any contrast function $ \rho $, $
f\in\bH_d(\beta,L,M) $, $ \delta>0 $ and $ n\in\bN^* $ such that 
$$
\frac{\varkappa_\delta}{\sqrt{N_b}}\geq 2\sup_{h\in[h_{\min},h_{\max}]}\big(1\vee b_h\sqrt{nh^d}\big),
$$
we have
$$
\bP_f\left[{\bar G_\delta^h}\right]
\leq
N_b\varSigma\:\exp\left\{-\frac{nh^d\left(\varkappa_\delta/2
\sqrt{N_b}\right)^2}{8
K_\infty^2\:(1\vee\dot{\rho}_\infty^2)+\frac{4\varkappa_\delta}{3\sqrt{N_b}}
\:K_\infty\:(1\vee{\dot{\rho}_\infty})}\right\}.
$$ 
\end{lemma}
Proofs of those lemmas are given in Appendix.

\paragraph{Proof of  Proposition \ref{lemma1}.}
Definitions of $\hat\theta(h)$ and
$\theta=\theta(f,x_0,h)$
imply that for any $\varepsilon\geq
\frac{4N_b}{c\lambda}\big(1\vee b_h\sqrt{nh^d}\big) $
\begin{eqnarray*}
\mathbb{P}_f\left(\sqrt{nh^d}\big|\hat{f}^h({x_0})-f({x_0})\big|
\geq\varepsilon,\:G_\delta^h\right)&\leq&\mathbb{P}
_f\left(\sqrt{nh^d}\big|\hat\theta_{0,\ldots,0}
(h)-\theta_{0,\ldots,0}\big|\geq
\varepsilon,\:G_\delta^h\right)\nonumber\\
&\leq&\mathbb{P}_f\left(\sqrt{nh^d}\sqrt{N_b}
\:\big\|\hat\theta(h)-\theta\big\|_2\geq
\varepsilon,\:G_\delta^h\right),
\end{eqnarray*}
where $ \|\cdot\|_2 $ is the $ \ell_2 $-norm on $ \bR^{N_b}
$.
Under the event $ G_\delta^h $ we have $
\hat\theta(h)\in\cB(\theta,\delta) $ for the specific choice of $\delta$ given in Lemma \ref{criterion_inversible} and
depending on the contrast $\rho$. According to Lemma \ref{criterion_inversible} {\it (2)} we obtain that
\begin{eqnarray*}
&&\mathbb{P}_f\left(\sqrt{nh^d}\big|\hat{f}
^h({x_0})-f({x_0})\big|\geq\varepsilon,\:G_\delta^h\right)\\
&&\quad\leq\mathbb{P}_f\left(\sqrt{nh^d}{\frac{2\sqrt{N_b}}{
c\lambda}}
\left\|D_h\big(\hat\theta(h)\big)-D_h(\theta)\right\|_2\geq
\varepsilon\right).
\end{eqnarray*}
Using Lemma \ref{criterion_inversible} {\it (1)} and the definition of $ \hat\theta(h) $ in (\ref{minimization}),
reminding that $ D_h(\theta)= \tilde
D_h\big(\hat\theta(h)\big)=0 $ and using the well-known
inequality $
\|\cdot\|_2\leq\sqrt{N_b}\|\cdot\|_\infty $ (where
$\|\cdot\|_2$ and $\|\cdot\|_\infty$ are respectively $ \ell_2 $-norm and $\ell_\infty$-norm on $
\bR^{N_b} $), we get with the last inequality:
\begin{eqnarray*}
&&\mathbb{P}_f\left(\sqrt{nh^d}\big|\hat{f}
^h({x_0})-f({x_0})\big|\geq\varepsilon,
\:G_\delta^h\right)\nonumber\\
&&\quad\leq\mathbb{P}_f\left(\sqrt{nh^d}{\frac{2\sqrt{N_b}}{
c\lambda}}\left\|\tilde
D_h\big(\hat\theta(h)\big)-D_h\big(\hat\theta(h)\big)\right\|_2\geq
\varepsilon\right)\nonumber\\
&&\quad\leq\sum_{ 
p\in\cS_b}\mathbb{P}_f\left(\sqrt{nh^d}\sup_{t\in\Theta(M)}
\left|\tilde D_h^{  p}(t)-D_h^{  p}(t)\right|\geq
\frac{c\lambda\:\varepsilon}{{2N_b}}
\right).\\
\end{eqnarray*}
Applying Lemma \ref{large_deviation} with $
z=\frac{c\lambda\:\varepsilon}{{2N_b}}$ and the
last inequality, finally we obtain the assertion
of Proposition \ref{lemma1}
\begin{eqnarray*}
&&\mathbb{P}_f\left(\sqrt{nh^d}\big|\hat{f}
^h({x_0})-f({x_0})\big|\geq\varepsilon,\:G_\delta^h\right)\\
&&\quad\leq
N_b\varSigma\exp\left\{-\frac{\left(\frac{
c\lambda\:\varepsilon}{{2N_b}}-
\big(1\vee b_h\sqrt{nh^d}\big) \right)^2}{8
K_\infty^2\:(1\vee\dot{\rho}_\infty^2)+\frac{4\:K_\infty}{3{N_b}}
\:(1\vee{\dot{\rho}_\infty})\frac{c\lambda\:
\varepsilon}{\sqrt{nh^d}}}\right\}.
\end{eqnarray*}
 \epr

\subsection{Proof of Theorem \ref{minimax}}\label{section_proof_maximalrisk}

Before starting Proofs of the main results of this paper, let us
define auxiliary results.
The next proposition provides us with upper bound for the
risk of a
$ \rho $-LPA estimator. Put
\begin{eqnarray}\label{def_constante_upperbound}
\bar C_r&=&\left(\frac{ 4
N_b}{{c\lambda}}\right)^r
+N_b\varSigma\int_{\frac{ 4
N_b}{{c\lambda}}}^{\infty}rz^{r-1}
\nonumber\\&&\times\exp\left\{-\frac{\left(\frac{z}{{N_b}}
{\frac{c\lambda}{2}}-1\right)^2}{8
K_\infty^2(1\vee\dot{\rho}_\infty^2)+\frac{4\delta}{3N_b}{
c\lambda}\:K_\infty\:(1\vee{\dot{\rho}_\infty})}\right\}\:dz,
\quad r\geq1.\qquad
\end{eqnarray}

\begin{proposition}\label{2}
Let $ \rho $ be a contrast function. Then, for any $ n\in\bN^*$, $ h>n^{-1/d} $, $ x_0\in[0,1]^d $ and any
$ f $ such that
$ \|\theta\|_1<\infty $, we have
\begin{eqnarray*}
\bE_f\big|\hat{f}^h({x_0})-f({x_0})\big|^r\bI_{G_\delta^h}\leq
\bar C_r\:\big(1\vee b_h\sqrt{nh^d}\big)^r\:\big(nh^d\big)^{-r/2},\quad
r\geq1.
\end{eqnarray*}
\end{proposition}
The proof of Proposition \ref{2} is
deduced from Proposition \ref{lemma1} by integration.

\paragraph{Proof of Theorem \ref{minimax}}

 By definition of $ \bH_d(\beta,L,M) $, the approximation error (bias) $b_h $
as defined in (\ref{def_term_bias}) verified $ b_h\leq Ldh^\beta $ for any $ h>0 $. Moreover by definition of $\bar
h=(L^2n)^{-\frac{1}{2\beta+d}}$ in (\ref{hbar}), we
have that
$b_{\bar h}\:\sqrt{n\bar h^d}\leq d$ and $ \big(n\bar
h^d\big)^{-1/2}=L^{\frac{d}{2\beta+d}}\varphi_n(\beta).$
We get
\begin{equation}\label{eq_ind}
 \mathbb{E}_f\big|\hat{f}^{\bar
h}({x_0})-f({x_0})\big|^r=\mathbb{E}_f\big|\hat{f}^{\bar
h}({x_0})-f({x_0})\big|^r\bI_{G_\delta^{\bar h}}
+\mathbb{E}_f\big|\hat{f}^{\bar
h}({x_0})-f({x_0})\big|^r\bI_{\bar G_\delta^{\bar h}}.
\end{equation}

The right hand side is controlled by Lemma
\ref{Control_complement}. Indeed, we can use the
Cauchy-Schwarz inequality,
\begin{eqnarray}\label{eq_risk_complement}
&&\mathbb{E}_f\big|\hat{f}^{\bar
h}({x_0})-f({x_0})\big|^r\bI_{\bar G_\delta^{\bar
h}}\nonumber\\
&&\quad\leq\left(\mathbb{E}_f\big|\hat{f}^{\bar
h}({x_0})-f({x_0})\big|^{2r}\bP_f\left\{\bar G_\delta^{\bar
h}\right\}\right)^{1/2}\nonumber\\
&&\quad\leq(2M)^r\sqrt{N_b\varSigma}\:\exp\left\{-\frac{
n\bar h^d\left(\varkappa_\delta
/2\sqrt{N_b}\right)^2}{16 K_\infty^2\:(1\vee\dot{\rho}_\infty^2)+\frac{8\varkappa_\delta}{3\sqrt{N_b}
}\:K_\infty\:(1\vee{\dot{\rho}_\infty})}\right\}.\qquad\quad
\end{eqnarray}
The last inequality is obtained because $M$ is a upper bound of $ f $ and $ \hat f^{\bar h} $ (cf. Definition
\ref{def_holder} and (\ref{ensemble})).
Using Proposition \ref{2}, (\ref{eq_ind}) and
(\ref{eq_risk_complement}), we obtain
\begin{eqnarray*}
&&\mathbb{E}_f\big|\hat{f}^{\bar
h}({x_0})-f({x_0})\big|^r\nonumber\\
&&\quad\leq \bar C_r
d^r\:L^{\frac{rd}{2\beta+d}}\varphi_n^r(\beta)\\
&&\quad\quad+(2M)^r\sqrt{N_b\varSigma}\:\exp\left\{-\frac{
n\bar h^d\left(\varkappa_\delta
/2\sqrt{N_b}\right)^2}{16 K_\infty^2\:(1\vee\dot{\rho}_\infty^2)+\frac{8\varkappa_\delta}{3\sqrt{N_b}
}\:K_\infty\:(1\vee{\dot{\rho}_\infty})}\right\}.
\end{eqnarray*}
When $ n $ tends towards $ +\infty $, Theorem \ref{minimax}
is proved. \epr
\subsection{Proof of Theorem \ref{lepski}}\label{section_proof_lepski} We start the
proof with formulating some auxiliary results whose proofs
are given in
Appendix. Define
\begin{equation}\label{fenetre_oracle}
 h^*=\left[\frac{\varrho_n(\beta)}{L^2d^2\:n}\right]^{\frac{1}{2\beta+d}},
\end{equation}
where $ \varrho_n^2(\beta) $ is defined in (\ref{normalizationAdaptive}).
Let $\kappa$ be an integer defined as follows:
\begin{equation}\label{kappa}
2^{-\kappa}h_{\max}\leq h^*<2^{-\kappa+1}h_{\max}.
\end{equation}
For  any $ n $ large enough, we have $ h_{\min}\leq
h^*\leq h_{\max} $.
\begin{lemma}\label{controlproba}
 For any $f\in\bH_d(\beta,L,M)$, any $ n $ large enough and any   $k\geq \kappa+1$
\begin{eqnarray*}
\bP_f\big(\hat{k}=k,\:G_\delta^{h_k}\big)\leq
J\:2^{-2(k-1)rd},
\end{eqnarray*}
where $ J=N_b\varSigma\big(1+(1-2^{-2rd})^{-1}\big) $.
\end{lemma}

\paragraph{Proof of Theorem \ref{lepski}.} This proof is
based on the scheme due to \cite{Lepski_Mammen_Spokoiny97}.
The  definition of $h^*$ (\ref{fenetre_oracle}) and
$\kappa$ (\ref{kappa}) implies that for any $n$ large enough
\begin{equation}\label{Control_produit}
\left(1\vee b_{h_k}\sqrt{nh_k^d}\right) \leq2\sqrt{1+k\ln2},\quad\forall
k\geq\kappa.
\end{equation}
Using Proposition \ref{2}, the last inequality yields
\begin{equation}\label{control_riskadaptive}
\bE_f\big|\hat{f}^{(k)}({x_0})-f({x_0})\big|^r\bI_{G_\delta^{h_k}}
\leq 
\bar C_r \:S_n^r(k),\quad\forall k\geq\kappa.
\end{equation}
To get this result we have applied Proposition \ref{2} with
$h=h_k$ and (\ref{Control_produit}).
We also have
\begin{eqnarray}\label{decomposition}
&&\mathbb{E}_f\big|\hat{f}^{(\hat
k)}({x_0})-f({x_0})\big|^r\nonumber\\
&&\quad=\mathbb{E}_f\big|\hat{f}^{(\hat
k)}({x_0})-f({x_0})\big|^r\bI_{\hat k\leq\kappa,G_\delta^{h_{\hat
k}}}+\mathbb{E}_f\big|\estimatorb{f}^{(\hat
k)}({x_0})-f({x_0})\big|^r\bI_{\hat k>\kappa,G_\delta^{h_{\hat
k}}}\nonumber\\
&&\qquad+\mathbb{E}_f\big|\estimatorb{f}^{(\hat
k)}({x_0})-f({x_0})\big|^r\bI_{\bar G_\delta^{h_{\hat
k}}}\nonumber\\
&&\quad:=R_1(f)+R_2(f)+R_3(f).
\end{eqnarray}
First we control $R_1$. By convexity of $ |\cdot|^r,\:r\geq1 $ and with the triangular inequality, we have
$$
\big|\hat{f}^{(\hat
k)}({x_0})-f({x_0})\big|^r\leq2^{r-1}\big|\hat{f}^{(\hat
k)}({x_0})-\hat{f}^{(\kappa)}({x_0})\big|^r+2^{r-1}\big|\estimatorb{
f}^{(\kappa)}({x_0})-f({x_0})\big|^r.
$$
The definition of $\hat k$ in (\ref{indexe adaptive}) yields
$$
\big|\hat{f}^{(\hat
k)}({x_0})-\hat{f}^{(\kappa)}({x_0})\big|^r\bI_{\hat
k\leq\kappa,G_\delta^{h_{\hat k}}}\leq C^rS_n^r(\kappa),
$$
where the constant $ C $ is defined in (\ref{variance}).
In view of
(\ref{control_riskadaptive}), the definitions of $ h_\kappa $
lead to
\begin{eqnarray*}
&&\bE_f\big|\hat{f}^{(\kappa)}({x_0})-f({x_0})\big|^r\bI_{\hat
k\leq\kappa,G_\delta^{h_{\hat k}}}\leq\bar
C_r\:S_n^r(\kappa),
\end{eqnarray*}
where $ \bar C_r $ is defined in (\ref{def_constante_upperbound}).
Noting that the right hand side of the obtained inequality is
independent of $f$ and taking into account the definition of
$\kappa$ and $h^*$ we obtain
\begin{eqnarray}\label{r200}
\limsup_{n\rightarrow\infty}\sup_{f\in\bH_d(\beta,L,A,M)}
\phi_n^{-r}(\beta)R_1(f)<\infty.
\end{eqnarray}
Now, let us bounded from above $R_2$. Applying
Cauchy-Schwartz
inequality, in view of Lemma \ref{controlproba} we have for $n$ large enough
\begin{eqnarray*}\label{r2}
R_2(f)&=&\sum_{k=\kappa}^{\mathbf{k}_n}\mathbb{E}
_f\big|\estimatorb{f}^{(k)}({x_0})-f({x_0})\big|^r\bI_{G_\delta^{h_{
k}}}\nonumber\\
&\leq&\sum_{k>\kappa}\left(\mathbb{E}_f\big|\estimatorb{f}^{
(k)}({x_0})-f({x_0})\big|^{2r}\bI_{G_\delta^{h_{
k}}}\right)^{1/2}\sqrt{\mathbb
P_f\big\{\hat k=k,G_\delta^{h_{
k}}\big\}}\nonumber\\
&=&\sqrt{J}\sum_{k>\kappa}\left(\mathbb{E}_f\big|\estimatorb
{f}^{(k)}({x_0})-f({x_0})\big|^{2r}\bI_{G_\delta^{h_{
k}}}\right)^{1/2}
2^{-(k-1)rd}.
\end{eqnarray*}
We obtain from
(\ref{control_riskadaptive}) and the last inequality
\begin{eqnarray*}
R_2(f)\leq
\frac{\bar{C}_{2r}^{r}\:2^{rd}\sqrt{J}
}{\big(nh_{\max}^d\big)^{r/2}}\:\sum_{s\geq
0}(1+s\ln2)^{r/2}2^{-srd}.
\end{eqnarray*}
 It remains to note that the right hand side of the last inequality
is  independent of $f$.
Thus, we have
\begin{eqnarray}
\label{r300}
&&\limsup_{n\to\infty}\sup_{f\in\bH_d(\beta,L,A,M)}\phi^{-r}
_n(\beta)R_2(f)<\infty.
\end{eqnarray}
It remains to bound $ R_3(f) $. By definition, note that  $ |\hat{f}^{(\hat k)}({x_0})|\leq M $, this allows us to state
that $\big|\hat{f}^{(\hat
k)}({x_0})-f({x_0})\big|\leq 2M$. Finally we obtain
$$
R_3(f)\leq
2^rM^r\:\bP_f\big\{\bar G_\delta^{h_{\hat k}}\big\}.
$$
Since $ nh_{\min}^d=\big(\ln n\big)^{2d} $, then
\begin{eqnarray}
\label{r400}
\limsup_{n\to\infty}\sup_{f\in\bH_d(\beta,L,A,M)}\phi^{-r}
_n(\beta)R_3(f)<\infty,
\end{eqnarray} follows now from  Lemma
\ref{Control_complement}. Theorem \ref{lepski} is proved
from (\ref{decomposition}), (\ref{r200}), (\ref{r300}) and (\ref{r400}).
 \epr

\section{Appendix}\label{appendixA}
\subsection{Proof of Lemma \ref{criterion_inversible}}
{\it 1.} By definition of the Jacobian matrix $ J_D(\cdot) $ in (\ref{Jacobian_matrix}), we can write
for any $ p, q\in\cS_b $ 
\begin{eqnarray*}
\left[J_D\big(\tilde\theta\big)\right]_{p,q}=
\frac{1}{n}\sum_{i=1}^n\int_{[-1/2,1/2]^d}x^{  p+ 
q}\:K(x)\int_\bR
\rho''\big(z-f_{\tilde\theta-\theta}
(y+hx)\big)\:g_i(z)dz\:dx.
\end{eqnarray*}
Applying this formula when $\tilde\theta=\theta$, the term $f_{\tilde\theta-\theta}$ vanishes, so:
\begin{eqnarray*}
 J_D\big(\theta\big)=
\frac{1}{n}\sum_{i=1}^n\int_\bR
\rho''(z)\:g_i(z)dz\:\int_{[-1/2,1/2]^d}U(x)\:U^\top
(x)\:K(x)dx.
\end{eqnarray*}
where $ U(\cdot) $ is defined in (\ref{objet_polynome}).
Since $ (g_i)_i\in\cG_{\rho}^{(c)} $, the definition of $ c $ in (\ref{def_constant_c}) implies that $
\frac{1}{n}\sum_{i=1}^n\int_\bR
\rho''(z)\:g_i(z)\:dz\geq c>0 $.

Now we show that
$ J_D(\theta) $ is a strictly
positive definite matrix, indeed for any $
\tau\in\bR^{N_b}\backslash0 $ 
\begin{eqnarray}\label{eq_lemma1_eigen_positive}
\tau^\top J_D(\theta)\tau&=&\frac{1}{n}\sum_{i=1}^n\int_\bR
\rho''(z)\:g_i(z)dz\:\tau^\top  
\int_{\big[-\frac{1}{2},\frac{1}{2}\big]^d}U(x)\:U^\top
(x)\:K(x)\:dx\:\tau\nonumber\\
&=&\frac{1}{n}\sum_{i=1}^n\int_\bR
\rho''(z)\:g_i(z)dz\:\int_{\big[
-\frac
{1}{2},\frac{1}{2}\big]
^d}\left[\tau^\top
U(x)\right]^2K(x)\:dx\nonumber\\
&\geq& c\int_{\big[
-\frac
{1}{2},\frac{1}{2}\big]
^d}\left[\tau^\top
U(x)\right]^2K(x)\:dx>0.
\end{eqnarray}

Let us show that for any $ h>n^{-1/d} $, $ \theta$ is the unique solution of $ D_h(\cdot)=(0,\ldots,0) $. By definition
in
(\ref{Deterministic_criterion}), $ D_h $ can be written as 
\begin{eqnarray}\label{Expression_Deterministic_criterion}
 D_h^{  p}(t)=-\frac{1}{n}\sum_{i=1}^n\int_{[-0.5,0.5]^d}
x^{ 
p}K(x)\int_\bR
\rho'\big(z-f_{t-\theta}(y+hx)\big)\:g_i(z)\;dz\;dx.\quad
\end{eqnarray}

Moreover, we have that 
$$
D_h(t)=(0,\ldots,0)\Longrightarrow\sum_{p\in\cS_b} \big(t_p-\theta_p\big)D_h^p(t)=0.
$$
Denote $ u(\cdot)= f_{t-\theta}(y+h\cdot)$. 
Since for any $i$, $ g_i $ is monotonically on $\bR_+$ and symmetric, then 
\begin{equation}\label{eq_lemma1_diff_densities_positive}
\forall
x\in[-0.5,0.5]^d,\quad\inf_{z>0}\inf_{i=1,\ldots,n}g_i\big(z-|u(x)|\big)-g_i\big(z+|u(x)|\big)\geq 0.
\end{equation}
Since $ (g_i)_i $ are symmetric, $ K $ is positive and $ \rho' $ is odd and positive on $ \bR_+^* $, the last inequality
and (\ref{Expression_Deterministic_criterion})
imply
\begin{eqnarray}\label{eq_lemma3_diff}
&&\int_{[-0.5,0.5]^d}
u(x)K(x)\int_\bR
\rho'\big(z-u(x)\big)\:\frac{\sum_{i=1}^ng_i(z)}{n}\;dz\;dx=0\nonumber\\
&\Leftrightarrow&\int_{[-0.5,0.5]^d}
u(x)K(x)\int_\bR
\rho'(z)\:\sum_{i=1}^ng_i\big(z+u(x)\big)\;dz\;dx=0\nonumber\\
&\Leftrightarrow&\int_{[-0.5,0.5]^d}
u(x)K(x)\int_0^{\infty}
\rho'(z)\:\sum_{i=1}^ng_i\big(z-u(x)\big)-g_i\big(z+u(x)\big)\;dz\;dx=0\quad\nonumber\\
&\Leftrightarrow&\int_{[-0.5,0.5]^d}
|u(x)|K(x)\int_0^{\infty}
\rho'(z)\:\sum_{i=1}^ng_i\big(z-|u(x)|\big)-g_i\big(z+|u(x)|\big)\;dz\;dx=0\quad\nonumber\\
&\Leftrightarrow&\forall
x\in[-0.5,0.5]^d,\:z>0,\quad\sum_{i=1}^ng_i\big(z-|u(x)|\big)-g_i\big(z+|u(x)|\big)=0\quad
\end{eqnarray}
Assume that there exists $ x_0\in [-0.5,0.5]^d$ such that $ u(x_0)\neq0 $. In particular for any $i$, since $
g_i $ is
monotonically on $\bR_+$, there exists $ z_{i,x_0}>0 $ such that
$$
g_i(z_{i,x_0}-|u(x_0)|)-g_i\big(z_{i,x_0}+|u(x_0)|\big)>0.
$$
That leads a contradiction in view of (\ref{eq_lemma1_diff_densities_positive}) and (\ref{eq_lemma3_diff}), thus for any
$x\in[-0.5,0.5]^d$, we have $ u(x)=0 $.
By definition of $ u(\cdot) $, we get
$$
\forall
x\in[-0.5,0.5]^d,\:h>n^{-1/d},\quad |f_{t-\theta}(y+hx)|=0\Longrightarrow t=\theta.
$$
Then, $ \theta $ is the unique solution of $ D_h(\cdot)=0 $.
%
%
%
%
%
%
\bigskip

{\it 2.} 
Let $ |||\cdot|||_2 $  be the euclidian matrix norm, $
\lambda_{\max}(A) $ the spectral ray of the matrix $A$ and
$ \lambda_0(A) $ the smallest eigenvalue of $A$.
According to Lemma \ref{criterion_inversible} {\it (1)}, there exists a radius $ \delta>0 $,
which only depends of $\rho$ such that
\begin{equation}\label{def_radius}
\inf_{\tilde\theta\in\cB(\theta,\delta)}\lambda_0\left(
J_D\big(\tilde\theta\big)
\right)\geq\lambda_0\left(J_D\big(\theta\big)\right)/2>0.
\end{equation}
This assertion can be explained as follows. In view of Assumption \ref{properties_rho} ({\it 3}),
$\lambda_0\left(J_D(\cdot)\right)$ is a continuous function. 
So, there exists a radius $ \delta>0 $ expected
such that (\ref{def_radius}) is true.

According to the local inverse function theorem, we can deduced that for any $ \tilde\theta\in\cB(\theta,\delta), $
\begin{equation}\label{eq_proof_lemme1_2}
 |||J_{D^{-1}}(\tilde\theta)|||_2=|||J_D^{-1}
(\tilde\theta)|||_2={\lambda_{\max}\big(J_D^{-1}
(\tilde\theta)\big)}=1/{\lambda_0\big(J_D
(\tilde\theta)\big)},
\end{equation}
By definition of $ \cG_{\rho}^{(c)} $, we have for any $ (g_i)_i\in\cG_{\rho}^{(c)} $, $ c=c(\rho,(g_i)_i) $ and
according to
(\ref{eq_lemma1_eigen_positive}) $ \lambda=\lambda_0\big(J_D
(\tilde\theta)\big)>0 $.
The smallest eigenvalue of $ J_D(\theta) $ is bigger than $
c\lambda>0 $. Indeed we have
\begin{eqnarray*}
 J_D\big(\theta\big)=
\frac{1}{n}\sum_{i=1}^n\int_\bR
\rho''(z)\:g_i(z)dz\int_{[-1/2,
1/2]^d}U(x)\:U^\top
(x)\:K(x)\:dx,
\end{eqnarray*}
and
\begin{eqnarray*}
\lambda_0\big(J_D(\theta)\big)&=&\frac{1}{n}\sum_{i=1}
^n\int_\bR
\rho''(z)\:g_i(z)dz\:\:\lambda_0\left(\int_{[-1/2,
1/2]^d}U(x)\:U^\top
(x)\:K(x)\:dx
\right)\\
&\geq& c\lambda.
\end{eqnarray*}

By definition of $ \delta $ in (\ref{def_radius}),
using (\ref{eq_proof_lemme1_2}) and the last inequality, we
have for any $
\tilde\theta\in\cB(\theta,\delta) $
\begin{equation}\label{eq_lemma1_borneJacobian}
  |||J_{D^{-1}}(\tilde\theta)|||_2\leq
{\frac{2}{c\lambda}}.
\end{equation}
As $ D_h $ is differentiable and each partial
derivative is continuous (cf. Assumption \ref{properties_rho} {\it (3)}), 
we use the local inverse function
theorem and (\ref{eq_lemma1_borneJacobian}) which give for
any $ \tilde\theta\in\cB(\theta,\delta) $ the following
inequality
$$
\|\tilde\theta-\theta\|_2=\left\|D_h^{-1}\circ
D_h\big(\tilde\theta\big)-D_h^{-1}\circ
D_h(\theta)\right\|_2\leq
{\frac{2}{c\lambda}
}\left\|
D_h\big(\tilde\theta\big)-
D_h(\theta)\right\|_2
$$\epr
\subsection{Proof of Lemma \ref{control_bias}} By definition
of $\cE_h^{  p}$ and $ D_h^{  p} $
in
(\ref{Deterministic_criterion}), we have for any $ t\in\Theta(M) $
\begin{eqnarray*}
\big|\cE_h^{  p}(t)-D_h^{  p}(t)\big|&\leq&\frac{1}{nh^d}\sum_{i=1}^n\int_{[0,1]^d}
\left|\frac{x-{x_0}}{h}\right|^{ 
p}\:K\left(\frac{x-{x_0}}{h}\right)\\
&&\times\int_\bR
\left|\rho'\big(z+f(x)-f_t(x)\big)-\rho'\big(z-f_{t-\theta}
(x)\big)\right|\:g_i(z)dz\:dx.
\end{eqnarray*}
Since $ \rho' $ is $ 1 $-Lipschitz (cf. Assumption
\ref{properties_rho} ({\it 2})) and $ \int K=1 $, then with the last
inequality, it yields
$$
\forall h>n^{-1/d},\quad\max_{p
\in\cS_b}\sup_{t\in\Theta(M)}\big|\cE_h^{  p}(t)-D_h^{  p}(t)\big|\leq b_h.
$$
\epr

\subsection{Proof of Lemma \ref{large_deviation}}

\paragraph{Bernstein's Inequality.}
To prove this lemma, we use the following well-known Bernstein's inequality which can be found in \cite{Massart07}
(Section 2.2.3, Proposition 2.9).
Let $ \cX_1,\ldots,\cX_n $ be independent square integrable random variables such that for some nonnegative constant $
\cX_\infty$, $\cX_i\leq\cX_\infty $ almost surely for all $ i=1,\cdots,n $. Then for any positive $ \epsilon $, we have
\begin{equation}\label{inequality_Bernstein}
 \bP\left(\sum_{i=1}^n\big(\cX_i-\bE\cX_i\big)\geq\epsilon\right)\leq\exp\left\{-\frac{\epsilon^2}{2\sum_{i=1}
^n\bE\cX_i^2
+2\cX_\infty\epsilon/3 } \right\},
\end{equation}
where $\bE=\bE^{n}$ is the mathematical expectation with
respect to the probability law $ \bP $ of $ \cX_1,\ldots,\cX_n $. The latter inequality is so-called {\it Bernstein's
inequality}.
\paragraph{Proof of Lemma \ref{large_deviation}.}

We have for any $   p\in\cS_b $
$$
\sup_{t\in\Theta(M)}\left|\tilde D_h^{  p}(t)-D_h^{ 
p}(t)\right|
\leq\sup_{t\in\Theta(M)}\left|\tilde D_h^{  p}(t)-\cE_h^{ 
p}(t)\right|+\sup_{t\in\Theta(M)}\big|\cE_h^{  p}(t)-D_h^{  p}(t)\big|.
$$
In view of Lemma \ref{control_bias}, we get 
\begin{equation}\label{eq1:lemma3}
\sup_{t\in\Theta(M)}\left|\tilde D_h^{  p}(t)-D_h^{ 
p}(t)\right|
\leq\sup_{t\in\Theta(M)}\left|\tilde D_h^{  p}(t)-\cE_h^{ 
p}(t)\right|+b_h.
\end{equation}
Set $L(\cdot)=\tilde D_h^{  p}(\cdot)-\cE_h^{  p}(\cdot)$.
To establish the assertion of the lemma, we use a chaining
argument on $
L(\cdot) $. Remember that $ \Theta(M) $ is a compact of $
\bR^{N_b} $ with $ \ell_1 $-norm. Let $ t_0\in\Theta(M) $ be
fixed and for any $ l\in\bN^* $
put $\Gamma_l$ a $
10^{-l}$-net on $ \Theta(M) $. We introduce
the following notations
\begin{eqnarray*}
 u_0(t)=t_0,\quad
u_l(t)=\arg\inf_{u\in\Gamma_l}\|u-t\|_1,\quad l\in\bN^*.
\end{eqnarray*}
Since $ \rho' $ is continuous, $ L(\cdot) $ is stochastically continuous which allows us to use 
the following chaining argument
\begin{eqnarray}\label{decompostion_chaining}
 L(t)=L(t_0)+\sum_{l=1}^\infty
L\big(u_l(t)\big)-L\big(u_{l-1}(t)\big),\quad \forall
t\in\Theta(M).
\end{eqnarray}
Using (\ref{eq1:lemma3}) and (\ref{decompostion_chaining}),
we obtain
\begin{eqnarray}\label{decompostion_chaining_eq2}
&&\bP_f\left(\sqrt{nh^d}\sup_{t\in\Theta(M)}\left|\tilde
D_h^{  p}(t)-D_h^{  p}(t)\right|\geq z\right)\nonumber\\
&&\quad\leq\bP_f\left(\sup_{t\in\Theta(M)}
\left|L(t)\right|\geq
\frac{z}{\sqrt{nh^d}}-b_h\right)\nonumber\\
&&\quad\leq\bP_f\left(\left|L(t_0)\right|+\sup_{
t\in\Theta(M)}\sum_{l=1}^\infty\left|
L\big(u_l(t)\big)-L\big(u_{l-1}(t)\big)\right|
\geq \frac{z}{\sqrt{nh^d}}-b_h\right).\qquad\quad
\end{eqnarray}
We can control the second term as follows.
$$
\sup_{t\in\Theta(M)}\sum_{l=1}^\infty\left|
L\big(u_l(t)\big)-L\big(u_{l-1}(t)\big)\right|
\leq\sum_{l=1}^\infty\sup_{\stackrel{u,
v\in\Gamma_l\times\Gamma_{l-1}}{\|u-v\|_1\leq 10^{-l}}}\left|
L(u)-L(v)\right|,
$$
where $ \Gamma_0=\{t_0\} $. Using
(\ref{decompostion_chaining_eq2}) and last inequality, we
get
\begin{eqnarray}\label{decomposition_chaining_eq3}
&&\bP_f\left(\sqrt{nh^d}\sup_{t\in\Theta(M)}
\left|L(t)\right|\geq z- b_h\:\sqrt{nh^d}\right)\nonumber\\
&&\quad\leq\bP_f\left(\sqrt{nh^d}\left|L(t_0)\right|\geq
z/2- b_h\:\sqrt{nh^d}/2\right)\nonumber\\
&&\qquad+\bP_f\left(\sqrt{nh^d}\sum_{l=1}^\infty\sup_{
\stackrel{u,v\in\Gamma_l\times\Gamma_{l-1}}{\|u-v\|_1\leq
10^{-l}}}\big| L(u)-L(v)\big|\geq
z/2- b_h\:\sqrt{nh^d}/2\right).\quad\qquad
\end{eqnarray}
By Definition of $
\tilde D_h^{  p} $ in (\ref{Derivate_rho_criterion}), we can write:
$$
\tilde D_h^{  p}(t)=\frac{1}{{nh^d}}\sum_{i=1}^n
\:\rho'\big(Y_i-f_t(X_i)\big)\:\left(\frac{X_i-{x_0}}{h}\right)^{
p}\:K\left(\frac{X_i-{x_0}}{h}\right).
$$
We define the function $
\cW_t(x,z)=\frac{1}{\sqrt{nh^d}}
\:\rho'\big(z+f(x)-f_t(x)\big)\:\left(\frac{x-{x_0}}{h}\right)^{
p}\:K\left(\frac{x-{x_0}}{h}\right) $ for all $ x\in[0,1]^d $ and $z\in\bR$.
Since for any $i$, $ Y_i=f(X_i)+\xi_i $, the process $ \sqrt{nh^d}\:L(\cdot) $ can be written as
an empirical process (sum of independent, zero-mean and bounded random
variables).
\begin{equation}\label{def_process_empiric}
\sqrt{nh^d}\:L(t)=\sum_{i=1}^n
\cW_t\big(X_i,\xi_i\big)-\bE_f
\cW_t\big(X_i,\xi_i\big),\quad t\in\Theta(M)
\end{equation}
At a fixed point $ t_0 $, we can use classical exponential
inequalities for empirical process. 
By definition of $ \cW_t(.,.) $ above, we have
\begin{equation}\label{property_process_empiric}
\sum_{i=1}^n\bE_f\cW_t^2(X_i,\xi_i)\leq \dot{\rho}_\infty^2\: K_\infty^2,\qquad \|\cW_t(.,.)\|_\infty\leq
{\dot{\rho}_\infty}\:
K_\infty/\sqrt{nh^d},
\end{equation}
where $ \|\cdot\|_\infty $ is the sup-norm. 

For the control of the first probability of
(\ref{decomposition_chaining_eq3}), we use the Bernstein's
inequality (\ref{inequality_Bernstein}), then
\begin{eqnarray}\label{decomposition_chaining_part1}
&&\bP_f\left(\sqrt{nh^d}\left|\tilde D_h^{  p}(t_0)-\cE_h^{ 
p}(t_0)\right|\geq
\frac{z}{2}-\frac{b_h\:\sqrt{nh^d}}{2}\right)\nonumber\\
&&\quad\leq2\exp\left\{-\frac{\left(z-
b_h\:\sqrt{nh^d}\right)^2}{8\dot{\rho}_\infty^2\:K_\infty^2+\frac{4{\dot{\rho}_\infty}\:K_\infty}{3\sqrt{nh^d}}\:(z-
b_h\:\sqrt{nh^d})}\right\}.
\end{eqnarray}
The second probability can be bounded as follows:
\begin{eqnarray}\label{decomposition_chaining_part3}
&&\bP_f\left(\sqrt{nh^d}\sum_{l=1}^\infty\sup_{\stackrel{u,
v\in\Gamma_l\times\Gamma_{l-1}}{\|u-v\|_1\leq 10^{-l}}}\big|
L(u)-L(v)\big|\geq
\frac{z}{2}-\frac{b_h\:\sqrt{nh^d}}{2}\right)\qquad\nonumber\\
 &&\leq\bP_f\left(\sqrt{nh^d}\sum_{l=1}^\infty\frac{1}{
l^2}\:\sup_{l\geq1}l^2\sup_{\stackrel{u,
v\in\Gamma_l\times\Gamma_{l-1}}{\|u-v\|_1\leq
10^{-l}}}\big|
L(u)-L(v)\big|\geq
\frac{z}{2}-\frac{b_h\:\sqrt{nh^d}}{2}\right)\nonumber\\
&&\leq\sum_{l=1}^\infty\sum_{\stackrel{u,
v\in\Gamma_l\times\Gamma_{l-1}}{\|u-v\|_1\leq
10^{-l}}}\bP_f\left(\sqrt{nh^d}\:\frac{\pi^2}{6}
l^2\:\big|L(u)-L(v)\big|\geq\frac{z}{2}-\frac{b_h\:\sqrt{
nh^d}}{2}\right).\qquad
\end{eqnarray}
In view of (\ref{def_process_empiric}), we notice that
$$
\sqrt{nh^d}\big[L(u)-L(v)\big]= \sum_{i=1}^n
\cW_u\big(X_i,\xi_i\big)-\cW_v\big(X_i,\xi_i\big)-\bE_f
\left[\cW_u\big(X_i,\xi_i\big)-\cW_v\big(X_i,
\xi_i\big)\right],
$$
then we have a sum of independent zero-mean random variables with finite variance and bounded.
Since $ \rho' $ is assumed Lipschitz, we have the following assertions.
\begin{eqnarray*}\label{property_process_empiric2}
\sum_{i=1}^n\bE_f\left[\cW_u(X_i,\xi_i)-\cW_v(X_i,
\xi_i)\right]^2&\leq&K_\infty^2 \|u-v\|_1^2,\nonumber\\
 \|\cW_u(.,.)-\cW_v(.,.)\|_\infty&\leq&K_\infty
\|u-v\|_1/\sqrt{nh^d},
\end{eqnarray*}
Using (\ref{decomposition_chaining_part3}), the Bernstein's inequality (\ref{inequality_Bernstein}) and  the last three
inequalities, we obtain
\begin{eqnarray}\label{decomposition_chaining_part2}
&&\bP_f\left(\sqrt{nh^d}\sum_{l=1}^\infty\sup_{\stackrel{u,
v\in\Gamma_l\times\Gamma_{l-1}}{\|u-v\|_1\leq 10^{-l}}}\big|
L(u)-L(v)\big|\geq
\frac{z}{2}-\frac{b_h\:\sqrt{nh^d}}{2}\right)\qquad\nonumber\\
&&\leq2\sum_{l=1}^\infty\sum_{\stackrel{u,
v\in\Gamma_l\times\Gamma_{l-1}}{\|u-v\|_1\leq
10^{-l}}}\exp\left\{-\frac{36\|u-v\|_1^{-1}}{\pi^4\:l^4}\right.\nonumber\\
&&\qquad\qquad\qquad\qquad\qquad\times\left.\frac
{\left(z-b_h\:\sqrt{nh^d}\right)^2}{8K_\infty^2\:\|u-v\|_1+\frac{4K_\infty}{3\sqrt{nh^d}}\:\left(z-b_h\:\sqrt{nh^d}
\right)}
\right\} \nonumber\\
 &&\leq2\sum_{l=1}^\infty\#(\Gamma_l)\:\#(\Gamma_{l-1}
)\:\exp\left\{-\frac{36\:10^l}{\pi^4\:l^4}\frac{
\left(z-b_h\:\sqrt{nh^d}\right)^2}{
8K_\infty^2+\frac{4K_\infty}{3\sqrt{nh^d}}\left(z-b_h\:\sqrt{nh^d}\right)}
\right\},\qquad\quad
\end{eqnarray}
where $ \#(\Gamma_l) $ is the cardinal of $ \Gamma_l $. Moreover, we notice that $ \#(\Gamma_l)\leq d10^{l} $.
Recall that $ z\geq 2\big(1\vee b_h\:\sqrt{nh^d}\big) $ and we notice that $
\displaystyle\min_{l\in\bN^*}\frac{18\:10^l}{\pi^4\:l^4}>1 
$. The last assertions allows us to write that
\begin{eqnarray*}
 &&\exp\left\{-\frac{36\:10^l}{\pi^4\:l^4}\frac{
\left(z-b_h\:\sqrt{nh^d}\right)^2}{
8K_\infty^2+\frac{4K_\infty}{3\sqrt{nh^d}}\left(z-b_h\:\sqrt{nh^d}\right)}
\right\}\\
&&\quad\leq \exp\left\{-\frac{18\:10^l}{\pi^4\:l^4}\frac{(8K_\infty)^{-1}}{K_\infty+1/3}
\right\}\times\exp\left\{-\frac{
\left(z-b_h\:\sqrt{nh^d}\right)^2}{
8K_\infty^2+\frac{4K_\infty}{3\sqrt{nh^d}}\left(z-b_h\:\sqrt{nh^d}\right)}
\right\}.
\end{eqnarray*}
Using
(\ref{decompostion_chaining_eq2}),
(\ref{decomposition_chaining_eq3}),
(\ref{decomposition_chaining_part1}), (\ref{decomposition_chaining_part2}) and the last inequality, we have for any $
p\in\cS_b $
\begin{eqnarray*}
&&\bP_f\left(\sqrt{nh^d}\sup_{t\in\Theta(M)}\left|\tilde
D_h^{  p}(t)-\cE_h^{  p}(t)\right|\geq z\right)\\
&&\quad\leq\varSigma\exp\left\{-\frac{\left(z-
b_h\:\sqrt{nh^d}\right)^2}{4K_\infty^2\:(1\vee\dot{\rho}_\infty^2)+\frac{4
K_\infty}{3\sqrt{nh^d}}\:(1\vee{\dot{\rho}_\infty})z}\right\},
\end{eqnarray*}
where $ \varSigma $ is defined in (\ref{def_constant_sigma}). This concludes the proof of Lemma
\ref{large_deviation}.
\epr
\subsection{Proof of Lemma \ref{Control_complement}}
Remember that the event $ {\bar G_\delta^h} $ can be
written as 
${\bar G_\delta^h}=\big\{\hat\theta(h)\notin\cB(\theta,
\delta)\big\}$
and $ \hat\theta(h) $ and $ \theta $ are respectively the
solutions of equations $ \tilde D_h(\cdot)=0
$ and $
D_h(\cdot)=0 $. Moreover $ \theta $ is the unique solution of $ D_h(\cdot)=0 $, then we can notice the following
inclusion
\begin{equation}\label{inclusion_outside}
 \big\{\hat\theta(h)\notin\cB(\theta,\delta)\big\}
\subseteq\left\{\sup_{t\in\Theta(M)\backslash\cB(\theta,
\delta)}\|\tilde
D_h(t)-D_h(t)\|_2\geq\varkappa_\delta\right\},
\end{equation}
where $ \displaystyle\varkappa_\delta=\inf_{
t\in\Theta(M)\backslash\cB(\theta,\delta)}\|D_h(t)\|_2/2 $. The latter inclusion can be interpreted as follows. In view
of Lemma \ref{criterion_inversible} {\it (1)}, $ \theta $ is the unique solution of $ D_h(\cdot)=0 $ thus $ D_h(\cdot) $
is not null on $\Theta(M)\backslash\cB(\theta,\delta)$. Moreover, $ D_h(\cdot) $ does not depend of $ n $, then $
\varkappa_\delta $ is positive and does not depend on $ n $. The event $ 
\big\{\hat\theta(h)\notin\cB(\theta,\delta)\big\} $ implies that there exists $
\tilde\theta\in\Theta(M)\backslash\cB(\theta,\delta) $ such that $ \tilde D_h\big(\tilde\theta\big)=0 $, then on a
neighborhood of $ \tilde\theta $, $ D_h(\cdot) $ and $ \tilde D_h(\cdot) $ are not closed. So, there exists $
\bar\theta\in\Theta(M)\backslash\cB(\theta,\delta) $ such that 
$$
\|\tilde
D_h\big(\bar\theta\big)-D_h\big(\bar\theta\big)\|_2\geq\varkappa_\delta.
$$
Then, the latter inequality implies (\ref{inclusion_outside}) by passing to the supremum.

Applying the inclusion (\ref{inclusion_outside}), we
obtain
$$
\bP_f\left({\bar G_\delta^h}\right)
\leq\sum_{ 
p\in\cS_b}\bP_f\left(\sqrt{nh^d}\sup_{
t\in\Theta(M)\backslash\cB(\theta,\delta)}\big|\tilde D_h^{ 
p}(t)-D_h^{ 
p}(t)\big|>\frac{\sqrt{nh^d}\varkappa_\delta}{\sqrt{N_b}}\right)
$$
Assumptions on $ n,h $ in Lemma \ref{Control_complement} allow us to show that 
$\sqrt{nh^d}\varkappa_\delta/\sqrt{N_b} \geq 2\big(1\vee\:b_h\:\sqrt{nh^d}\big).$
Using Lemma \ref{large_deviation} with $ z=\sqrt{nh^d}\varkappa_\delta/\sqrt{N_b} $, we have 
$$
\bP_f\left({\bar G_\delta^h}\right)\leq
N_b\varSigma\:\exp\left\{-\frac{
nh^d\left(\varkappa_\delta/2\sqrt{N_b}\right)^2}{8
K_\infty^2\:(1\vee\dot{\rho}_\infty^2)+\frac{4\varkappa_\delta}{3\sqrt{N_b}}
\:K_\infty\:(1\vee{\dot{\rho}_\infty})}\right\}.
$$
The lemma is proved.\epr
\subsection{Proof of Lemma \ref{controlproba}}
Note that by definition of
$\hat k$ in
(\ref{indexe adaptive})
$$
\forall k\geq \kappa+1,\quad\big\{\hat{k}=k\big\}=\cup_{l\geq
k}\left\{\big|\estimatorb{f}^{(k-1)}({x_0})-\hat{f}^{(l)}({x_0})\big|> C\: S_n(l)\right\}.
$$
Note that $S_n(l)$ is monotonically increasing in $l$ and,
therefore,
\begin{eqnarray*}
\big\{\hat{k}=k\big\}&\subseteq &\left\{\big|\estimatorb{f}^{(k-1)}({x_0})-f({x_0})\big|> 2^{-1}C\: S_n(k-1)\right\}
\\*[2mm]
& &\cup\left[ \cup_{l\geq k}\left\{\big|\hat{f}^{(l)}({x_0})-f({x_0})\big|>
2^{-1}C\: S_n(l)\right\}\right].
\end{eqnarray*}
We come to the following
inequality: for any $k\geq \kappa+1$
\begin{eqnarray}
\label{eq1:proof-lemma2}
\bP\left(\hat{k}=k,G_\delta^{h_k}\right)&\leq&\bP\left\{\big|\hat{f}^{(k-1)}({x_0})-\hat{f}({x_0})\big|>
2^{-1}C\; S_n(k-1),G_\delta^{h_k}\right\} \nonumber\\*[2mm] & &+
\sum_{l\geq k}\bP\left\{\big|\hat{f}^{(l)}({x_0})-f({x_0})\big|> 2^{-1}C\;
S_n(l),G_\delta^{h_k}\right\}.\qquad
\end{eqnarray}
Notice that the definition of $S_n(l)$ yields
$N_{h_l}\;S_n(l)=
\big[1+l\ln2\big]^{1/2}.$
Thus, applying Proposition \ref{lemma1} with
$\ve=C\big[1+l\ln2\big]^{1/2}$ and $h=h_l$ and using the inequality (\ref{Control_produit}),
we obtain by definition of $ C $ in (\ref{variance}), for any $ l\geq k-1$ and $ n $ large enough
\begin{eqnarray}
\label{eq2:proof-lemma2} \bP\left\{\big|\hat{f}^{(l)}({x_0})-f({x_0})\big|>
2^{-1}C\; S_n(l)\right\}&\leq&N_b\varSigma\:2^{-2rd\:l}.
\end{eqnarray}
We obtain from
(\ref{eq1:proof-lemma2}) and  (\ref{eq2:proof-lemma2}) that $k\geq \kappa+1$
\begin{eqnarray*}
&&\bP\left(\hat{k}=k,G_\delta^{h_k}\right)\leq J2^{-2(k-1)rd},
\end{eqnarray*}
where $J=N_b\varSigma\big(1+(1-2^{-2rd})^{-1}\big)$. \epr

{\footnotesize
\bibliographystyle{plainnat}
\bibliography{reference2}

\begin{thebibliography}{40}
\providecommand{\natexlab}[1]{#1}
\providecommand{\url}[1]{\texttt{#1}}
\expandafter\ifx\csname urlstyle\endcsname\relax
  \providecommand{\doi}[1]{doi: #1}\else
  \providecommand{\doi}{doi: \begingroup \urlstyle{rm}\Url}\fi

\bibitem[Astola et~al.(2010)Astola, Egiazarian, Foi, and
  Katkovnik]{Katkovnik_Foi_Egiazarian_Astola10}
J.~Astola, K.~Egiazarian, A.~Foi, and V.~Katkovnik.
\newblock From local kernel to nonlocal multiple-model image denoising.
\newblock \emph{Int. J. Comput. Vision}, 86\penalty0 (1):\penalty0 1--32, 2010.

\bibitem[Barron et~al.(1999)Barron, Birg{\'e}, and
  Massart]{Barron_Birge_Massart99}
A.~Barron, L.~Birg{\'e}, and P.~Massart.
\newblock Risk bounds for model selection via penalization.
\newblock \emph{Probab. Theory Related Fields}, 113\penalty0 (3):\penalty0
  301--413, 1999.

\bibitem[Bousquet(2002)]{Bousquet02}
O.~Bousquet.
\newblock A {B}ennett concentration inequality and its application to suprema
  of empirical processes.
\newblock \emph{C. R. Math. Acad. Sci. Paris}, 334\penalty0 (6):\penalty0
  495--500, 2002.

\bibitem[Brown and Low(1996)]{Brow_Low96}
L.~Brown and M.~Low.
\newblock A constrained risk inequality with applications to nonparametric
  functional estimation.
\newblock \emph{Ann. Statist.}, 24\penalty0 (6):\penalty0 2524--2535, 1996.

\bibitem[Brown et~al.(2008)Brown, Cai, and Zhou]{Brown_Cai_Harrison08}
L.D. Brown, T.~Tony Cai, and H.H. Zhou.
\newblock Robust nonparametric estimation via wavelet median regression.
\newblock \emph{Ann. Statist.}, 36\penalty0 (5):\penalty0 2055--2084, 2008.

\bibitem[Chang and Guo(2005)]{Chang_Guo}
X.W. Chang and Y.~Guo.
\newblock Huber's m-estimation in relative gps positioning: computational
  aspects.
\newblock \emph{Journal of Geodesy.}, 2005.

\bibitem[Chichignoud(2011)]{Chichignoud10}
M.~Chichignoud.
\newblock Minimax and minimax adaptive estimation in multiplicative regression
  : locally bayesian approach.
\newblock \emph{Probab. Theory Related Fields}, 2011.
\newblock to appear.

\bibitem[Donoho et~al.(1995)Donoho, Johnstone, Kerkyacharian, and
  Picard]{Donoho_Johnstone_Kerkyacharian_Picard95}
D.~Donoho, I.~Johnstone, G.~Kerkyacharian, and D.~Picard.
\newblock Wavelet shrinkage: asymptopia?
\newblock \emph{J. Roy. Statist. Soc. Ser. B}, 57\penalty0 (2):\penalty0
  301--369, 1995.
\newblock With discussion and a reply by the authors.

\bibitem[Goldenshluger and Lepski(2008)]{Goldenshluger_Lepski08}
A.~Goldenshluger and O.V. Lepski.
\newblock Universal pointwise selection rule in multivariate function
  estimation.
\newblock \emph{Bernoulli}, 14\penalty0 (3):\penalty0 1150--1190, 2008.

\bibitem[Goldenshluger and Lepski(2009)]{Goldenshluger_Lepski09}
A.~Goldenshluger and O.V. Lepski.
\newblock Structural adaptation via lp-norm oracle inequalities.
\newblock \emph{Probab. Theory and Related Fields}, 143:\penalty0 41--71, 2009.

\bibitem[Goldenshluger and Nemirovski(1997)]{Goldenshluger_Nemirovski97}
A.~Goldenshluger and A.~Nemirovski.
\newblock On spatially adaptive estimation of nonparametric regression.
\newblock \emph{Math. Methods Statist.}, 6\penalty0 (2):\penalty0 135--170,
  1997.

\bibitem[Hall and Jones(1990)]{Hall_Jones89}
Peter Hall and M.~C. Jones.
\newblock Adaptive {$M$}-estimation in nonparametric regression.
\newblock \emph{Ann. Statist.}, 18\penalty0 (4):\penalty0 1712--1728, 1990.

\bibitem[H{\"a}rdle and Tsybakov(1988)]{Hardle_Tsybakov88}
W.~H{\"a}rdle and A.B. Tsybakov.
\newblock Robust nonparametric regression with simultaneous scale curve
  estimation.
\newblock \emph{Ann. Statist.}, 16\penalty0 (1):\penalty0 120--135, 1988.

\bibitem[H{\"a}rdle and Tsybakov(1992)]{Hardle_Tsybakov92}
W.~H{\"a}rdle and A.B. Tsybakov.
\newblock Robust locally adaptive nonparametric regression.
\newblock In \emph{Data analysis and statistical inference}, pages 127--144.
  Eul, Bergisch Gladbach, 1992.

\bibitem[Huber and Ronchetti(2009)]{Huber_Ronchetti09}
P.~Huber and E.~Ronchetti.
\newblock \emph{Robust statistics}.
\newblock Wiley Series in Probability and Statistics. John Wiley \& Sons Inc.,
  Hoboken, NJ, second edition, 2009.

\bibitem[Huber(1964)]{Huber64}
P.J. Huber.
\newblock Robust estimation of a location parameter.
\newblock \emph{Ann. Math. Statist.}, 35:\penalty0 73--101, 1964.

\bibitem[Juditsky(1997)]{Judistsky97}
A.~Juditsky.
\newblock Wavelet estimators: adapting to unknown smoothness.
\newblock \emph{Math. Methods Statist.}, 6\penalty0 (1):\penalty0 1--25, 1997.

\bibitem[Juditsky et~al.(2009)Juditsky, Lepski, and
  Tsybakov]{Juditsky_Lepski_Tsybakov09}
A.B. Juditsky, O.V. Lepski, and A.B. Tsybakov.
\newblock Nonparametric estimation of composite functions.
\newblock \emph{Ann. Statist.}, 37\penalty0 (3):\penalty0 1360--1404, 2009.

\bibitem[Katkovnik(1985)]{Katkovnik85}
V.~Katkovnik.
\newblock \emph{Nonparametric identification and data smoothing}.
\newblock ``Nauka'', Moscow (in Russian), 1985.
\newblock The method of local approximation.

\bibitem[Kerkyacharian et~al.(2001)Kerkyacharian, Lepski, and
  Picard]{Kerkyacharian_Lepski_Picard01}
G.~Kerkyacharian, O.V. Lepski, and D.~Picard.
\newblock Non linear estimation in anisotropic multi-index denoising.
\newblock \emph{Probab. Theory and Related Fields}, 121:\penalty0 137--170,
  2001.

\bibitem[Klutchnikoff(2005)]{Klutchnikoff05}
N.~Klutchnikoff.
\newblock \emph{On the adaptive estimation of anisotropic functions}.
\newblock PhD thesis, Aix-Masrseille 1, 2005.

\bibitem[Lepski et~al.(1997)Lepski, Mammen, and
  Spokoiny]{Lepski_Mammen_Spokoiny97}
O.~V. Lepski, E.~Mammen, and V.~G. Spokoiny.
\newblock Optimal spatial adaptation to inhomogeneous smoothness: an approach
  based on kernel estimates with variable bandwidth selectors.
\newblock \emph{Ann. Statist.}, 25\penalty0 (3):\penalty0 929--947, 1997.

\bibitem[Lepski(1990)]{Lepski90}
O.V. Lepski.
\newblock On a problem of adaptive estimation in gaussian white noise.
\newblock \emph{Theory of Probability and its Applications}, 35\penalty0
  (3):\penalty0 454--466, 1990.

\bibitem[Lepski(1991)]{Lepski91}
O.V. Lepski.
\newblock Asymptotically minimax adaptive estimation i. upper bounds. optimally
  adaptive estimates.
\newblock \emph{Theory Probab. Appl.}, 36:\penalty0 682--697, 1991.

\bibitem[Lepski and Levit(1999)]{Lepski_Levit99}
O.V. Lepski and B.Y. Levit.
\newblock Adaptive nonparametric estimation of smooth multivariate functions.
\newblock \emph{Mathematicals methods of statistics}, 1999.

\bibitem[Lepski and Spokoiny(1997)]{Lepski_Spokoiny97}
O.V. Lepski and V.G. Spokoiny.
\newblock Optimal pointwise adaptive methods in nonparametric estimation.
\newblock \emph{Annals of statistics}, 25\penalty0 (6):\penalty0 2512--2546,
  1997.

\bibitem[Massart(2000)]{Massart00}
P.~Massart.
\newblock Some applications of concentration inequalities to statistics.
\newblock \emph{Probability Theory}, volume spécial dédié à Michel
  Talagrand\penalty0 (2):\penalty0 245--303, 2000.

\bibitem[Massart(2007)]{Massart07}
P.~Massart.
\newblock \emph{Concentration inequalities and model selection}, volume 1896 of
  \emph{Lecture Notes in Mathematics}.
\newblock Springer, Berlin, 2007.
\newblock Lectures from the 33rd Summer School on Probability Theory held in
  Saint-Flour, July 6--23, 2003, With a foreword by Jean Picard.

\bibitem[Petrus(1999)]{Petrus99}
P.~Petrus.
\newblock Robust huber adaptive filter.
\newblock \emph{IEEE Transactions on Signal Processing.}, 47:\penalty0
  1129--1133, 1999.

\bibitem[Reiss et~al.(2011)Reiss, Rozenholc, and
  Cuenod]{Reiss_Rozenholc_Cuenod09}
M.~Reiss, Y.~Rozenholc, and C.~Cuenod.
\newblock Pointwise adaptive estimation for robust and quantile regression.
\newblock 2011.
\newblock Source: Arxiv.

\bibitem[Rousseeuw and Leroy(1987)]{Rousseuw_Leroy87}
P.~Rousseeuw and A.~Leroy.
\newblock \emph{Robust regression and outlier detection}.
\newblock Wiley Series in Probability and Mathematical Statistics: Applied
  Probability and Statistics. John Wiley \& Sons Inc., New York, 1987.

\bibitem[Talagrand(1996a)]{Talagrand96a}
M.~Talagrand.
\newblock New concentration inequalities in product spaces.
\newblock \emph{Inventiones Mathematicae}, 126:\penalty0 505–563, 1996a.

\bibitem[Talagrand(1996b)]{Talagrand96b}
M.~Talagrand.
\newblock A new look at independence.
\newblock \emph{Annals of Probability}, 24:\penalty0 1–34, 1996b.

\bibitem[Tsybakov(1982a)]{Tsybakov82a}
A.~B. Tsybakov.
\newblock Nonparametric signal estimation when there is incomplete information
  on the noise distribution.
\newblock \emph{Problems of Information Transmission}, 18\penalty0
  (2):\penalty0 116--130, 1982a.

\bibitem[Tsybakov(1982b)]{Tsybakov82b}
A.~B. Tsybakov.
\newblock Robust estimates of a function.
\newblock \emph{Problems of Information Transmission}, 18\penalty0
  (3):\penalty0 39--52, 1982b.

\bibitem[Tsybakov(1983)]{Tsybakov83}
A.~B. Tsybakov.
\newblock Convergence of nonparametric robust algorithms of reconstruction of
  functions.
\newblock \emph{Automation and Remote Control}, \penalty0 (12):\penalty0
  66--76, 1983.

\bibitem[Tsybakov(1986)]{Tsybakov86}
A.~B. Tsybakov.
\newblock Robust reconstruction of functions by a local approximation method.
\newblock \emph{Problems of Information Transmission}, 22\penalty0
  (2):\penalty0 69--84, 1986.

\bibitem[Tsybakov(2008)]{Tsybakov08}
A.~B. Tsybakov.
\newblock \emph{Introduction to Nonparametric Estimation}.
\newblock Springer Publishing Company, Incorporated, 2008.

\bibitem[Tsybakov(1998)]{Tsybakov98}
A.B. Tsybakov.
\newblock Pointwise and sup-norm sharp adaptive estimation of function on the
  sobolev classes.
\newblock \emph{Annals of statistics}, 26\penalty0 (6):\penalty0 2420--2469,
  1998.

\bibitem[Van~der Vaart and Wellner(1996)]{VanderVaart_Wellner96}
Aad~W. Van~der Vaart and Jon~A. Wellner.
\newblock \emph{Weak convergence and empirical processes}.
\newblock Springer Series in Statistics. Springer-Verlag, New York, 1996.
\newblock With applications to statistics.

\end{thebibliography}
}

\end{document}